\RequirePackage{tikz}
\documentclass[sn-mathphys]{sn-jnl}

\usepackage{subcaption}

\jyear{2021}%

\theoremstyle{thmstyleone}%
\theoremstyle{thmstyletwo}%

\theoremstyle{thmstylethree}%

\begin{document}

\title[A model for reversible electroporation]{A model for reversible electroporation to deliver drugs into diseased tissues}


\author*[1]{\fnm{Nilay} \sur{Mondal}}\email{nilay.mondal@iitg.ac.in}
\equalcont{These authors contributed equally to this work.}

\author[1]{\fnm{D. C.} \sur{Dalal}}\email{durga@iitg.ac.in}
\equalcont{These authors contributed equally to this work.}


\affil*[1]{\orgdiv{{Department of Mathematics}, \orgname{Indian Institute of Technology Guwahati}, \orgaddress{\street{Amingaon}, \city{Guwahati}, \postcode{781039}, \state{Assam}, \country{India}}}}

%


\abstract{Drug delivery through electroporation could be highly beneficial for the treatment of different types of diseased tissues within the human body. In this work, a  mathematical model of reversible tissue electroporation is presented for injecting drug into the diseased cells. The model emphasizes the tissue boundary where the drug is injected as a point source. Drug loss from the tissue boundaries through extracellular space is studied. Multiple pulses are applied to deliver a sufficient amount of drug into the targeted cells. The set of differential equations that model the physical circumstances are solved numerically. This model obtains a mass transfer coefficient in terms of pore fraction coefficient and drug permeability. It controls the drug transport from extracellular to intracellular space. The drug penetration throughout the tissue is captured for the application of different pulses. The boundary effects on drug concentration are highlighted in this study. The advocated model is able to perform homogeneous drug transport into the cells so that the affected tissue is treated completely. This model can be applied to optimize clinical experiments by avoiding the lengthy and costly in vivo and in vitro experiments.}




\maketitle

\section{Introduction}
In the current situation, one of the important challenges is to deliver a sufficient amount of drugs or genes into cells for bio-therapy and cancer chemotherapy. In mass transport, the cell membrane does not allow large size particles to enter into the intracellular space. Different physical approaches, such as, electroporation, micro-injection, laser, ultrasound \cite{Bolhassani2011} have been proposed by the researchers. Electroporation is a successful tool to increase mass transfer rate into the targeted cells. In electroporation, the cell membrane gets permeabilized, and nanometer-sized pores are formed due to the application of external electric pulses \cite{Kotnik2019,Miklavcic2020,Krassowska2007}. The transitory and permeabilized state of the cell membrane can be utilized to transport different molecules, such as drugs, ions, dyes, antibodies, oligonucleotides, RNA and DNA into the intracellular domain \cite{Miklavcic2012,Miklavcic2018SSept}. Electroporation has wide applications both \textit{in vitro} and \textit{in vivo}. Electroporation is also used for the delivery of drugs to malignant tumors and other diseased tissues \cite{Kotnik2019}. 

According to the application of pulse strength and pulse length, electroporation occurs in two types: reversible and irreversible. Reversible electroporation leads to temporary permeabilization of cell membranes with the survival of the treated cells \cite{Garcia2014}. Cells that are reversibly electroporated endure the treatment, and after the pulse period is over, the cell membrane starts to reseal \cite{Davalos2003}. Membrane resealing is a physiological process, and it is observed that the time duration for pore formation is in the order of microsecond. In contrast, membrane resealing happens in a minute time frame \cite{Pavlin,Pavlin2003}. Sometimes, electroporation guides cell death due to permanent permeabilization of the cell membrane when the applied electric field is high enough. This results in consequent loss of cell homeostasis, and the process is termed as irreversible electroporation \cite{Rubinsky2007,Davalos2005}. Irreversibly electroporated cells can not sustain the treatment, and pores do not reseal \cite{Jiang2015}.  

The cellular drug uptake depends on the electroporation parameters, such as pulse strength, pulse duration, and the number of pulses \cite{Corovic,Pucihar2011}. The mass transfer rate also depends on the tissue conductivity, nature of cells, and the characterization of the drug \cite{Pucihar2011}.  Drug transport rate increases with the increase of pulse duration and the number of pulses. The molecules can enter into the cells through the permeabilized membrane during the time break between the pulses \cite{Satkauskas}. Usually, two types of electroporation pulse are prevalent, short and long. From experimental studies, it is observed that the application of short duration (1$\mu$s-100 $\mu$s) high voltage electric field increases pore density within the cell membrane \cite{Weaver2003,Pavlin}. Nanometer-sized pores are created within individual lipids through the application of a short pulse  \cite{Becker2013}. Bulk electroporation is generally conducted to electroporate the whole tissue by applying some long duration (100 $\mu$s - 100 ms) low voltage pulses \cite{Miklavcic2018Feb,Granot}. It is reasonable to consider that increment of permeability, resealing of pores, and cell death seems to regulate the drug uptake in reversibly electroporated cells. Researchers endeavored to develop mathematical models on the basis of the aforementioned factors.

Granot and Rubinsky \cite{Granot} first attempted to develop a mathematical model for drug delivery in the tissue cells with reversible electroporation and proposed MTC in terms of pore density using the model \cite{Krassowska2007}. A detailed relation between electroporation and the increment in the permeability that affects cellular drug uptake is theoretically described in a dual-porosity model by Kalamiza et al. \cite{Kalamiza20141950}. The mass transport coefficient depending on membrane electropores is illustrated to solute transport into the tissue in their study.  Argus et al. \cite{Becker2017}  developed a model where reversible and irreversible electroporations were considered for mass transport into different types of electroporated cells. The tissue electrical conductivity depending on electric field is studied in their work.  Goldberg et al.  \cite{Goldberg2018} also used the pore creation model and proposed a multiphysics model for ion transport. In their model, they presented a mass transport equation using the Nernst-Planck equation for transporting different species into the cells.

The recent study \cite{Mondal2021} proposes a mathematical model of tissue electroporation for transporting drugs into cells. This paper discusses both reversible and irreversible electroporation with uniform electric fields and pores resealing. The delivery of drugs during and after electroporation are depicted separately. There is no consideration about  the effects of drug loss or drug elimination from the tissue  in the above mentioned models. The method of drug injection as a point source is also not studied so far. 

In the present study, a mathematical model for drug delivery is developed in which bulk electroporation on tissue is taken into account. Only reversible electroporation has been emphasized in this model to cure whole infected tissue by delivering drug as a medicine. On this account, some low voltage multiple pulses are applied to succeed a sufficient amount of drugs into the targeted cells. The nobility of this study is that  a point source drug in terms of Dirac-delta function is introduced to represent initial drug distribution. In addition, the effects of drug loss from the tissue boundary on drug transport are studied thoroughly.  In this study, suitable differential equations with appropriate initial and boundary conditions, which  govern the biological circumstances, are presented. The main objective of the current study is to analyze the drug transport phenomena in the tissue by taking care of the effects of drug loss, which is imposed as boundary conditions. 

\section{Problem Formulation}

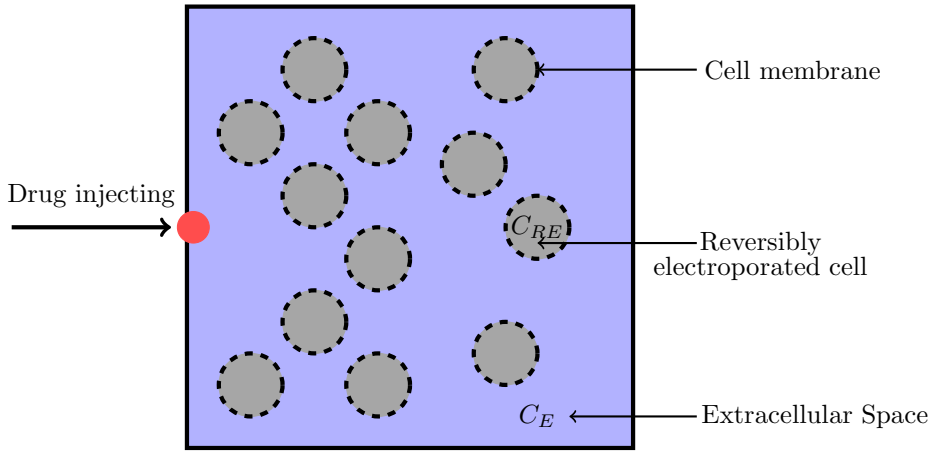
\begin{figure}[h!]
	\centering
	\begin{tikzpicture}[scale=0.42]
	
	\draw[fill,blue!30] (-6,-6) rectangle (8,8);
	\draw[ultra thick] (-6,-6) rectangle (8,8);
	
	\draw[fill,gray!70](0,0) circle (1cm);
	\draw[ultra thick, dashed] (0,0) circle (1cm);
	
	\draw[fill,gray!70] (-2,2) circle (1cm);
	\draw[ultra thick,dashed] (-2,2) circle (1cm);
	
	\draw[fill,gray!70](3,3) circle (1cm);
	\draw[ultra thick, dashed] (3,3) circle (1cm);
	
	\draw[fill,gray!70] (-4,4) circle (1cm);
	\draw[ultra thick,dashed] (-4,4) circle (1cm);

	\draw[fill,gray!70](4,-3) circle (1cm);
	\draw[ultra thick, dashed] (4,-3) circle (1cm);
	
	\draw[fill,gray!70] (4,6) circle (1cm);
	\draw[ultra thick,dashed] (4,6) circle (1cm);
	
	\draw[fill,gray!70](-4,-4) circle (1cm);
	\draw[ultra thick, dashed] (-4,-4) circle (1cm);
	
	\draw[fill,gray!70] (-2,-2) circle (1cm);
	\draw[ultra thick,dashed] (-2,-2) circle (1cm);
	
	
	\draw[fill,gray!70] (-2,6) circle (1cm);
	\draw[ultra thick,dashed] (-2,6) circle (1cm);
	
	\draw[fill,gray!70] (5,1) circle (1cm);
	\draw[ultra thick,dashed] (5,1) circle (1cm);
	
	\draw[fill,gray!70] (0,4) circle (1cm);
	\draw[ultra thick,dashed] (0,4) circle (1cm);
	
	\draw[fill,gray!70] (0,-4) circle (1cm);
	\draw[ultra thick,dashed] (0,-4) circle (1cm);
	
	\draw[thick, <-] (6,-5) -- (10,-5);
	\draw (13.7,-5) node[]{Extracellular Space};
	\draw (5,-5) node[]{$C_E$};
	
	\draw[thick, <-] (5,0.5) -- (10,0.5);
	\draw (12,0.5) node{Reversibly };
	\draw (12,-0.3) node{electroporated cell};
	\draw (5,1) node{$C_{RE}$};
	
	\draw[thick, <-] (5,6) -- (10,6);
	\draw (13,6) node{Cell membrane };
	
	\draw[ultra thick, ->] (-11.5,1) -- (-6.5,1);
	\draw (-9,2) node{Drug injecting};
	
	\draw[fill,red!70] (-5.8,1) circle (0.5cm);

	\end{tikzpicture}
	
	\caption{A schematic diagram of injecting drug  into a biological tissue.}
	\label{3tissue1}	
\end{figure}
\begin{figure}[h!]
	\centering
	\begin{tikzpicture}[scale=0.4]

	\draw[ultra thick] (-6,-6) rectangle (8,8);

	%
	\draw[ultra thick,->](-8,2) -- (-8,8.3);
	\draw[ultra thick,->](-8,0) -- (-8,-6.3);
	\draw (-8,1) node[]{\large L};
	\draw[ultra thick,->](6,2) -- (6,0);
	\draw (5,1) node[]{\Large E};	
	
	\draw[fill,black!100] (-6,8) rectangle (8,8.3);
	\draw[fill,black!100] (-6,-6) rectangle (8,-6.3);
	
	\draw (1,-8) node[]{\bf Electrode (-)};
	\draw (1,-5) node[]{\large $\phi=\phi_0$};
	
	\draw (1,10) node[]{\bf Electrode (+)};
	\draw (1,7) node[]{\large $\phi=\phi_L$};
	
	\draw[ultra thick,->](9,5) -- (9,8);
	\draw[ultra thick,->](9,5) -- (12,5);
	\draw (12.5,5) node[]{\large x};
	\draw (9,8.5) node[]{\large y};
	\draw (8.8,4.6) node[]{\large o};
	
	\draw (1,1) node[]{\large $\Omega$};

	\end{tikzpicture}
	
	\caption{A schematic representation of bulk electroporation on the square tissue.}
	\label{3electrode1}
\end{figure}
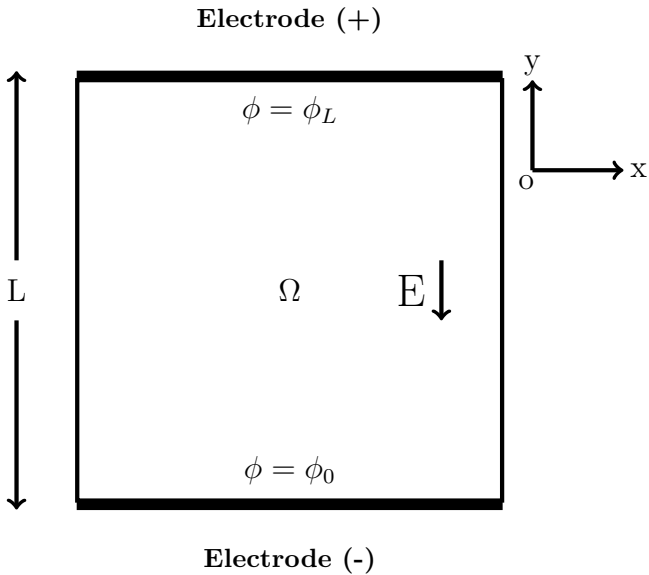
%
%
%
%
%
%
%
%
%
%

In this study, a biological tissue is considered and assumed to be cubical, having side length  $L$. Structurally, the tissue is split into two parts: extracellular space and intracellular space. Both the spaces are separated by the porous cell membranes, which control the mass transport phenomenon from extracellular to intracellular domains. A schematic diagram as shown in Fig. \ref{3tissue1} is provided to visualize the complete tissue structure as well as the drug injecting process into the tissue.  Two electrodes are set at the upper and lower boundaries of the tissue for electroporation. The positive electrode is assigned with the potential $\phi_L$, and the negative electrode is assigned with the potential $\phi_0$ (as shown in Fig. \ref{3electrode1}). Due to the application of electric pulses, the transmembrane potential is increased, and the cell membrane gets permeabilized. As a result, nanometer-sized pores are formed in the cell membrane and the pore fraction coefficient ($f_p$) increases. Repeated pulses are applied to the tissue over a specific time period with some time gap between the pulses. The pore resealing occurs in the rest time when a pulse is off. During this rest time, drug transport occurs under pore resealing effects. In the resealing period, the pore area decreases with time, and the drug transport rate gets slow.  Multiple pulses are applied after completing the rest time of the previous pulse to enhance the mass transfer rate. 

\section{Model Development and Mathematical Equations}
Because of the arrangement of two parallel electrodes, a uniform electric field and a non-uniform electric potential are generated. The potential ($\phi$) distribution inside the tissue is obtained by solving the equation \cite{Sel2005,Corovic,Bradley2016}
\begin{eqnarray}\label{3eq1}
\vec{\nabla}\cdot(\sigma \cdot \vec{\nabla} \phi)=0,
\end{eqnarray}
with boundary conditions: $\phi(y=0)=\phi_0$ and $\phi(y=L)=\phi_L$, where $\sigma$ is the tissue electrical conductivity.\\
The uniform electric field ($E$), through out the domain $\Omega$ is calculated by the magnitude of the gradient of the potential, mathematically  expressed as 
\begin{eqnarray}\label{3E}
E= \lvert \vec{\nabla} \phi \lvert
\end{eqnarray}
\subsection{Roles of tissue conductivity}
The creation of pores in the cell membrane is one of the effective outcomes of electroporation. In macroscopic view, these pores have the effect of increasing the electrical conductivity. Therefore, the electrical conductivity  increases as the electric field ($E$) increases. Previous studies have considered various functional dependencies of the electrical conductivity on the electric field magnitude \cite{Sel2005,Lackovic2010,Corovic2013}. From the  experimental as well as computational study \cite{Sel2005}, it is known that during electroporation the electrical conductivity changes with the applied electric field and the relation is defined as,
\begin{equation}\label{3eq3}
\sigma(E)=\left[\frac{\sigma_{\text{max}}-\sigma_{\text{min}}}{1+\gamma_1 \cdot \exp  \left(-\frac{E-a}{b}\right)}\right]+\sigma_{\text{min}} ,
\end{equation}
\begin{equation*}
a=\frac{E_{\text{rev}}+E_{\text{irrev}}}{2} \mbox{  and } b=\frac{E_{\text{irrev}}-E_{\text{rev}}}{\gamma_2},
\end{equation*}
where $E_{\text{rev}}$ and $E_{\text{irrev}}$ are respective threshold values for reversible and irreversible electroporation,
$\sigma_{\text{min}}$ and $\sigma_{\text{max}}$ are the minimum and the maximum tissue electrical conductivity respectively, $\gamma_1$ and $\gamma_2$ are the sigmoidal functional parameters.
\subsection{Mass transfer coefficient calculation}
The induced electric field is a cause for increasing transmembrane potential and as a result cell membrane get permeabilized.   
A symmetric sigmoid function, is developed in the model \cite{Becker2017} to evaluate pore fraction coefficient ($f_p$) of the surviving cells (reversibly electroporated cells) which is defined as,
\begin{align}\label{3eq4}
f_p=\left[1+\exp\left(\frac{E_{f}-E}{b_f}\right)\right]^{-1},
\end{align}
where $E_f$ and $b_f$ are the fitting parameters.\\
The mathematical formulation of the MTC ($\mu$) is given as follows,
\begin{equation}\label{3mu}
\mu(t)=\frac{Pf_p}{r_c}e^{-\frac{t}{\tau}},
\end{equation}
where $P$ is the drug permeability, $r_c$ is the cell radius and $\tau$ is the resealing constant.
\subsection{Drug transport phenomenon in the tissue}
The drug concentrations in extracellular space and in reversibly electroporated cells are obtained by the mass transport equation as
\begin{eqnarray}\label{3eq6}
\frac{\partial C_E}{\partial t}=\vec{\nabla}.\left(D \vec{\nabla} C_E\right) - \left(\frac{1-\varepsilon}{\varepsilon}\right)\mu (t)\times
\left(C_E - C_{RE}\right),
\end{eqnarray} 
\begin{equation}\label{3eq7}
\frac{\partial C_{RE}}{\partial t}=\mu (t)\times (C_E - C_{RE}),
\end{equation}
where  $C_E$ and $C_{RE}$ are the drug concentrations in extracellular space and reversible electroporated cells respectively; $D$ is the effective diffusion coefficient of the drug in the extracellular space; $\varepsilon$ is the porosity of the membrane pores. We solve the above equations considering initial and boundary conditions, which are as follows:
\begin{equation}\label{3eq8}
\left.\begin{split}
&C_E(x, y, 0)=C_1(x,y),\qquad C_{RE}(x, y, 0)=0,\\
&\left(\frac{\partial C_E}{\partial x}\right)_{x=0}=\beta C_E,\qquad \qquad
\left(\frac{\partial C_E}{\partial x}\right)_{x=L}=\beta C_E,\\
&\left(\frac{\partial C_E}{\partial y}\right)_{y=0}=\beta C_E,\qquad \qquad
\left(\frac{\partial C_E}{\partial y}\right)_{y=L}=\beta C_E,
\end{split}\qquad\right\}
\end{equation}
where $\beta$ represents the rate of mass loss from the tissue boundary which depends on the tissue's boundary properties. The initial drug that is injected through a particular point/location of the tissue boundary is mathematically defined as 
\begin{equation*}
C_1(x,y)=
\begin{cases}
n_d\delta(y-0.5),& x=0\\
\qquad 0,&\text{otherwise}
\end{cases}
\end{equation*}
where $n_d$ is the number of drug dose and $\delta(y)$ refers to the Dirac delta function. The delta function can be used in terms of normal distribution as
\begin{equation}\label{3eq9}
\delta_d(y)=\frac{1}{d\sqrt{\pi}} \exp^{\left(-\frac{y}{d}\right)^2}, \quad 0<d<1.
\end{equation}
\section{Method of Solution} 	  	
The governing equations \eqref{3eq1} - \eqref{3eq9} are solved  numerically. Equations \eqref{3eq1} - \eqref{3eq3} are solved simultaneously to find optimum potential and electric field throughout the tissue region. In numerical computations, the calculation is started on taking $\sigma=\sigma_{\text{min}}$ initially at t=0. A uniform electric field is induced in the tissue due to parallel electrode configuration at tissue boundaries. The field direction is from top to bottom, which is shown in the  Fig. \ref{3electrode1}. 
We can determine the pore fraction coefficient from the expression \eqref{3eq4} using the optimized electric field. The mass transfer coefficient is now calculated from the relation \eqref{3mu}.
The coupled equations \eqref{3eq6} - \eqref{3eq9} are solved numerically using  finite difference method. The  details of the parameters used in the model are shown in  Table \ref{3tab1}. 

\begin{table}[h!]
	\caption{Details of the parameters used in the simulation:}
	\vspace{0.2cm}
	\label{3tab1}
	\centering
	\begin{tabular}{llll}
		\hline
		Symbol & Value  & Description  & Source   \\
		\hline 
		&&\\
		$\sigma_{\text{min}}$ & 0.0 S m$^{-1}$ & Minimum electrical conductivity & \cite{Becker2017}  \\
		$\sigma_{\text{max}}$ & 0.241 S m$^{-1}$ & Maximum electrical conductivity & \cite{Becker2017}  \\
		$E_{\text{rev}}$ & 46 V mm$^{-1}$ & Threshold value for reversible electroporation & \cite{Becker2017}  \\
		$E_{\text{irrev}}$ & 70 V mm$^{-1}$ & Threshold value for irreversible electroporation & \cite{Becker2017}  \\
		$\gamma_1$  & 8 & Electrical conductivity parameter  & \cite{Becker2017}  \\
		$\gamma_2$  & 10 & Electrical conductivity parameter  & \cite{Becker2017}  \\
		$E_f$  & 65.8 V mm$^{-1}$ & Fitting parameter for MTC & \cite{Becker2017}  \\
		$b_f$  & 7.5 V mm$^{-1}$ & Fitting parameter for MTC & \cite{Becker2017}  \\
		$r_c$ & 50 $\mu$m & Cell radius & \cite{Krassowska2007}  \\
		$D$ & $10^{-3}$ mm$^2$ s$^{-1}$ & Effective diffusion coefficient & \\
		$\varepsilon$& $0.18$ & Porosity of the membrane pores & \cite{Kalamiza20141950} \\
		$P$ & $5\times 10^{-4}$ mm s$^{-1}$ & Permeability of drug & \cite{Granot} \\
		$E$ & 60 V mm$^{-1}$ & Electrical field & \\
		$L$ & $1$ mm  & Length of the rectangle & Fig. \ref{3electrode1} \\
		$\phi_L$ & 60 V  & Potential on positive electrode & Fig. \ref{3electrode1}\\
		$\phi_0$ & 0 V  & Potential on negative electrode & Fig. \ref{3electrode1}\\
		$PN$ & 10 & Pulse number& \\
		$t_{ep}$ & 1 ms & Pulse length (ON TIME) & \\
		$t_{M}$ & 100 s & Time for mass transfer (OFF TIME) & \\
		$n_{d}$ & 100 & Number of drug dose & \\
		$d$ & 0.1 & Constant in Dirac-delta function & Eq. \eqref{3eq9}\\
		$\beta$ & 0.1 mm$^{-1}$ & Rate of mass loss at the boundary & \\
		&&\\
		\hline
	\end{tabular}
\end{table}
In discretization, we have used forward time centered space (FTCS) scheme to the equations \eqref{3eq6} - \eqref{3eq9} and discretized equations are as follows,
\begin{align}\label{3eq10}
\left(C_E\right)_{i,j}^{n+1}=&a\left(C_E\right)_{i+1,j}^{n}+a\left(C_E\right)_{i-1,j}^{n}+b(t)\left(C_E\right)_{i,j}^{n} \\ \nonumber
&+c\left(C_E\right)_{i,j+1}^{n}+c\left(C_E\right)_{i,j-1}^{n}+d(t)\left(C_{RE}\right)_{i,j}^{n}, \quad 
\begin{array}{c}
i=1,2,\cdots, M_1\\
j=1,2,\cdots, M_2
\end{array}
\end{align}
\begin{align}\label{3eq11}
\left(C_{RE}\right)_{i,j}^{n+1}=\left(C_{RE}\right)_{i,j}^{n}+\mu(t)\Delta t\left[\left(C_E\right)_{i,j}^{n}-\left(C_{RE}\right)_{i,j}^{n}\right], \quad \begin{array}{c}
i=1,2,\cdots, M_1 \\
j=1,2,\cdots, M_2
\end{array}
\end{align}
where $a=\frac{D\Delta t}{(\Delta x)^2}$, 
$b(t)=1-\left[2D\left(\frac{1}{(\Delta x)^2}+\frac{1}{(\Delta y)^2}\right)+\frac{1-\varepsilon}{\varepsilon}\mu(t)\right]\Delta t$, 
$c=\frac{D\Delta t}{(\Delta y)^2}$, and $d(t)=\frac{1-\varepsilon}{\varepsilon}\mu(t)\Delta t$. Here, $\Delta t$ is the time step size, $\Delta x$ and $\Delta y$ are the step sizes for space, $M_1$ and $M_2$ are the numbers of grid along $x$-axis and $y$-axis respectively.\\
\textbf{\textit{Stability condition:}}
To solve the equations \eqref{3eq10} - \eqref{3eq11} using FTCS scheme, the following stability condition is used.
\begin{align*}
\Delta t<\frac{1}{2}\times \frac{(\Delta x)^2 (\Delta y)^2}{D[(\Delta x)^2+ (\Delta y)^2]}.
\end{align*}
Here, $M_1=101$, $M_2=101$, $\Delta x=0.01$, $\Delta y=0.01$ and $\Delta t=0.2$ are chosen for simulations.

\section{Results and Discussion}
\begin{figure}[h!]
	\centering
	\begin{subfigure}{.5\textwidth}
		\centering
		\includegraphics[width=\linewidth]{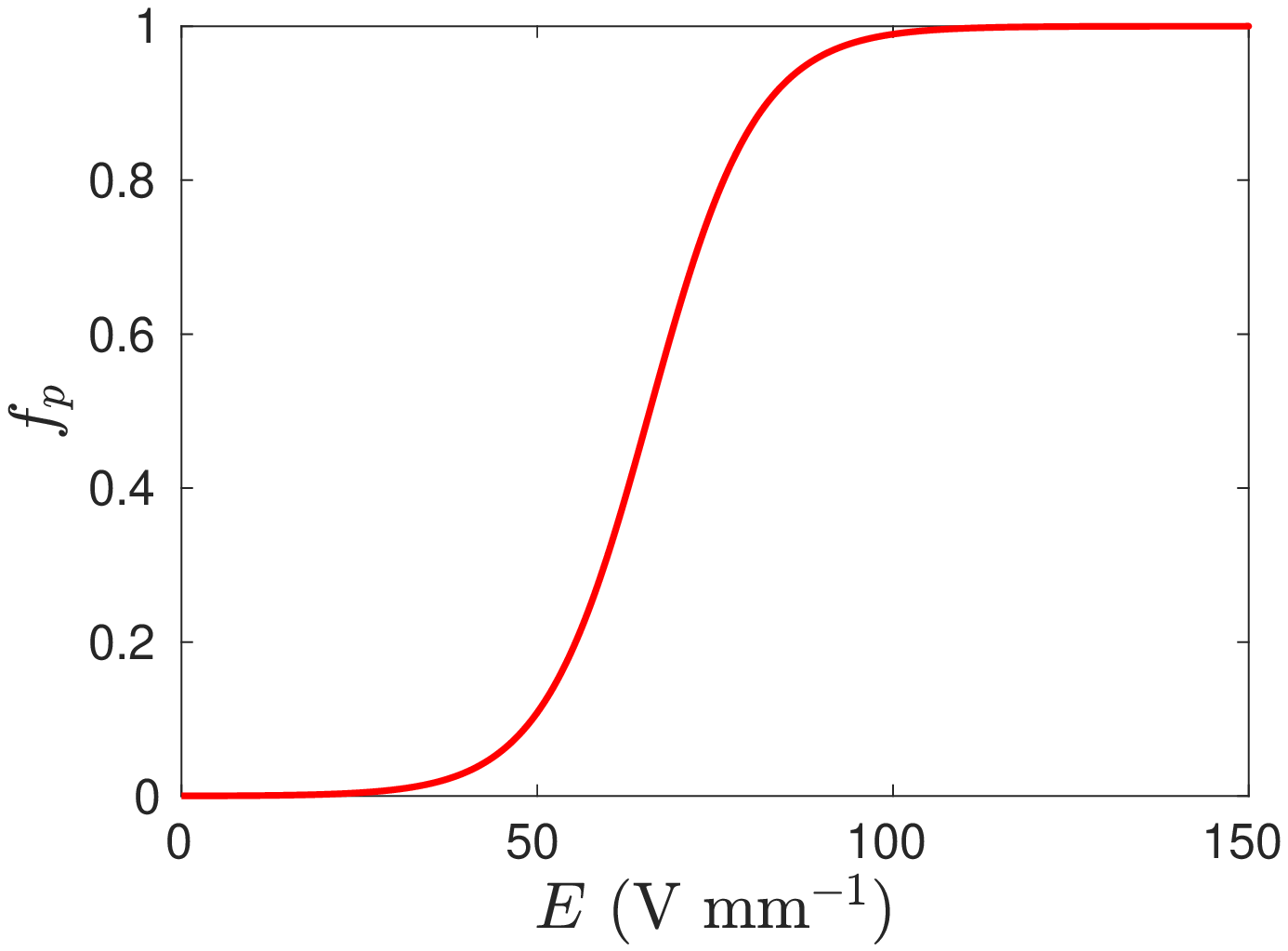}
		\caption{}
		\label{34sub1}
	\end{subfigure}%
	\begin{subfigure}{.5\textwidth}
		\centering
		\includegraphics[width=\linewidth]{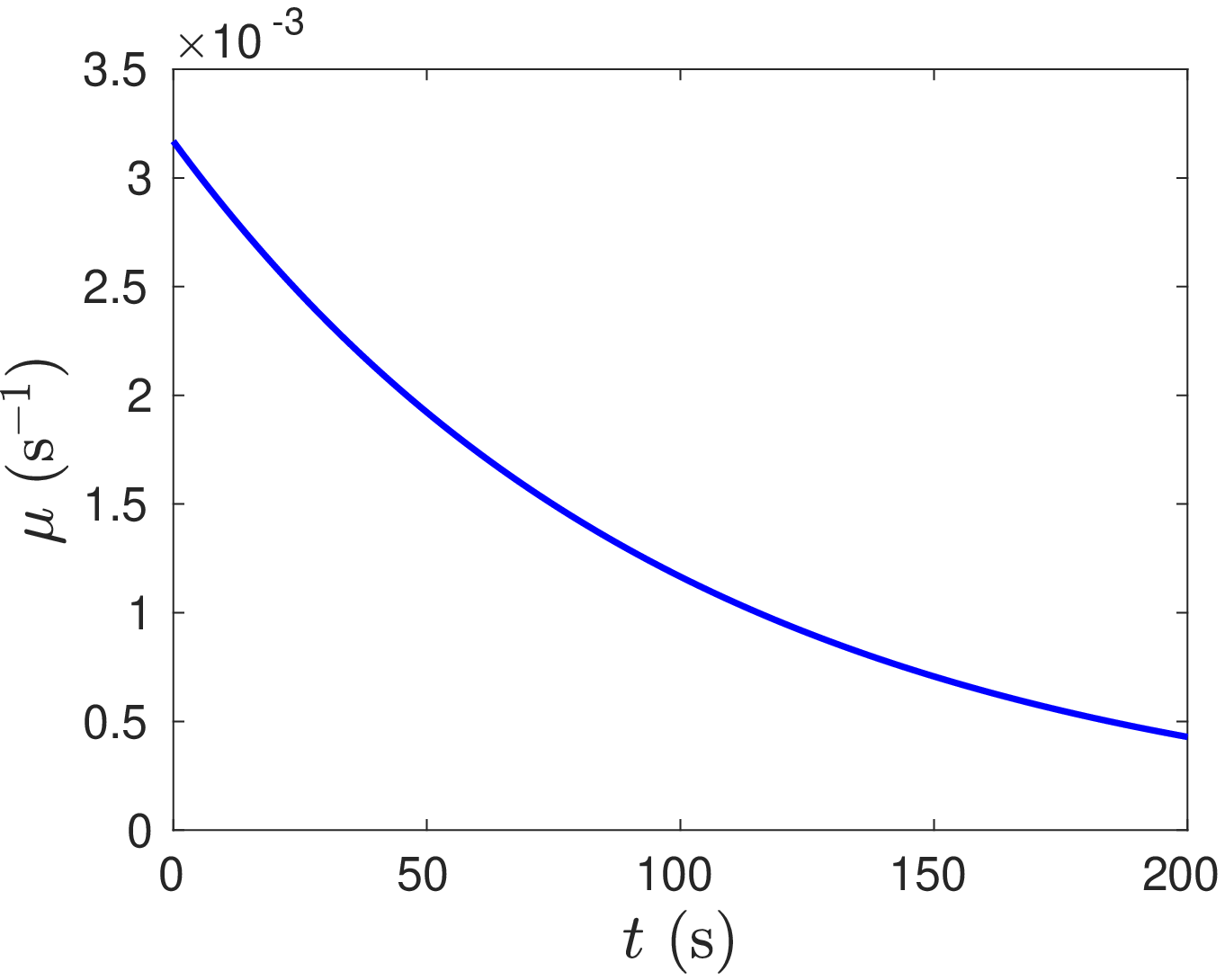}
		\caption{}
		\label{34sub2}
	\end{subfigure}
	\caption{Represents (a) the changes of pore fraction coefficient ($f_p$) with the applied electric field and (b) the changes of MTC with time. }
	\label{3fp_mu}
\end{figure}

In this study, a mathematical model and its numerical simulation are used to portray the drug transport phenomenon in diseased tissues.  A qualitative analysis is conducted to visualize how the drug is transported into the electroporated tissue using graphical representations given in the Figs. \ref{3fp_mu} - \ref{3CE_CRE_PN_change_beta_0}. The results are discussed through the plotted graphs for different parameters like, rate of drug loss ($\beta$), drug permeability ($P$) and pulse number ($PN$).

The drug transport from extracellular to intracellular space depends on the mass transfer rate through the cell membrane. This mass transfer rate depends on the pore fraction coefficient, which is an outcome of electroporation. The changes of $f_p$ with the electric field are shown by the graph in Fig. \ref{34sub1}. It is observed that the pore fraction coefficient suddenly increases when the applied field crosses a certain limit (50 V mm$^{-1}$). The maximum value ($=1$) of $f_p$ is attained when the field is greater than 100 V mm$^{-1}$. In this model, The strength of the electric field is considered to be  60 V mm$^{-1}$ to get a higher value of $f_p$ that enhances the mass transfer rate in the cell membrane. Also, this electric field does not cause irreversible electroporation as $E<E_{rev}$.  Fig. \ref{34sub2} shows that the MTC decreases with time due to the effect of membrane resealing after completion of a pulse. If MTC  decreases during mass transfer (OFF TIME), the drug uptake by the cells will be affected. Thus, in order to maintain a higher mass transfer rate, it is necessary to reopen the membrane pores by repeating pulses.

\subsection{Drug transport phenomenon in the tissue}
\begin{figure}[h!]
	\centering
	\begin{tikzpicture}[scale=0.35]
	
	\draw[ultra thick] (-6,-6) rectangle (8,8);
	
	\draw[ultra thick, ->] (-11,1) -- (-6,1);
	\draw (-8,2) node{Drug};
	\draw (-10.5,-0.5) node{$C_E=n_d\delta(y-0.5)$ };
	\draw (1,2) node{$C_E=0$};
	\draw (1,0) node{$C_{RE}=0$};
	\draw (1,9) node{$\frac{\partial C_E}{\partial y}=\beta C_E$ };
	\draw (1,-7) node{$\frac{\partial C_E}{\partial y}=\beta C_E$ };
	\draw (11,1) node{$\frac{\partial C_E}{\partial x}=\beta C_E$ };
	\end{tikzpicture}
	
	\caption{A schematic diagram of the model equations, which are presented in the Eqs. \eqref{3eq6}-\eqref{3eq9} for injecting drug  into a biological tissue.}
	\label{3model}	
\end{figure}
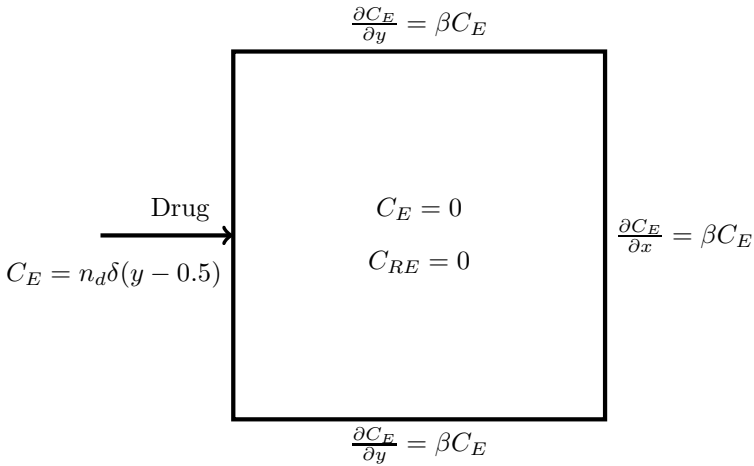

In this section, a complete overview about the drug transport from the point source to the tissue concerning the effects of drug losses from the boundaries of the tissue is discussed. The physical phenomenon of drug delivery with the prescribed boundary condition is highlighted through the picture given in the Fig. \ref{3model}.

\begin{figure}[h!]
	\centering
	\includegraphics[width=0.5\linewidth]{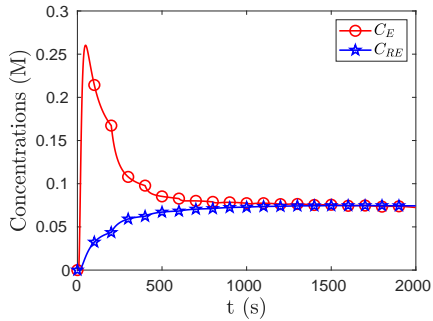}					
	\caption{Temporal changes of drug concentrations ($C_E$, $C_{RE}$). Here, concentrations are obtained for $\beta=0.1$ mm$^{-1}$ at the point (0.5, 0.5).}
	\label{3cons_x5}
\end{figure}
Fig. \ref{3cons_x5} represents the changes of drug concentration profiles in the ECS and ICS with time. In the figure, extracellular concentration initially increases sharply due to the diffusion of drug into the ECS from its source. Once the drugs start entering into the cells, the concentration decreases until it reaches an equilibrium state.  On the other hand, intracellular concentration increases with time due to continuous drug intake into the cells.  It is observed that both the profiles merge after a certain time due to drug saturation in the ECS and ICS.  
\subsubsection{Comparison}
\begin{figure}[h!]
	\centering
	\begin{tabular}{cc}
		\includegraphics[width=.5\linewidth]{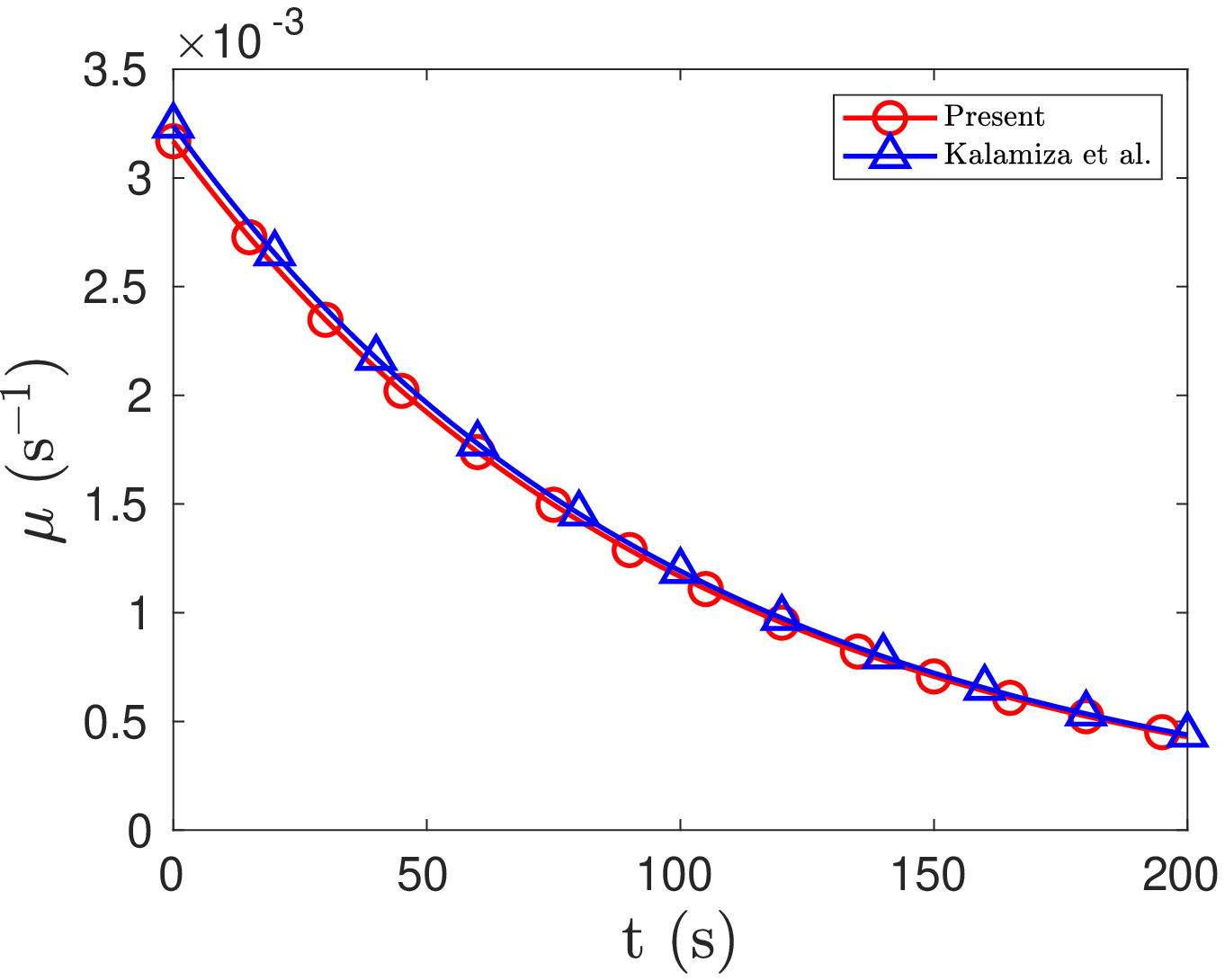}
		&
		\includegraphics[width=.5\linewidth]{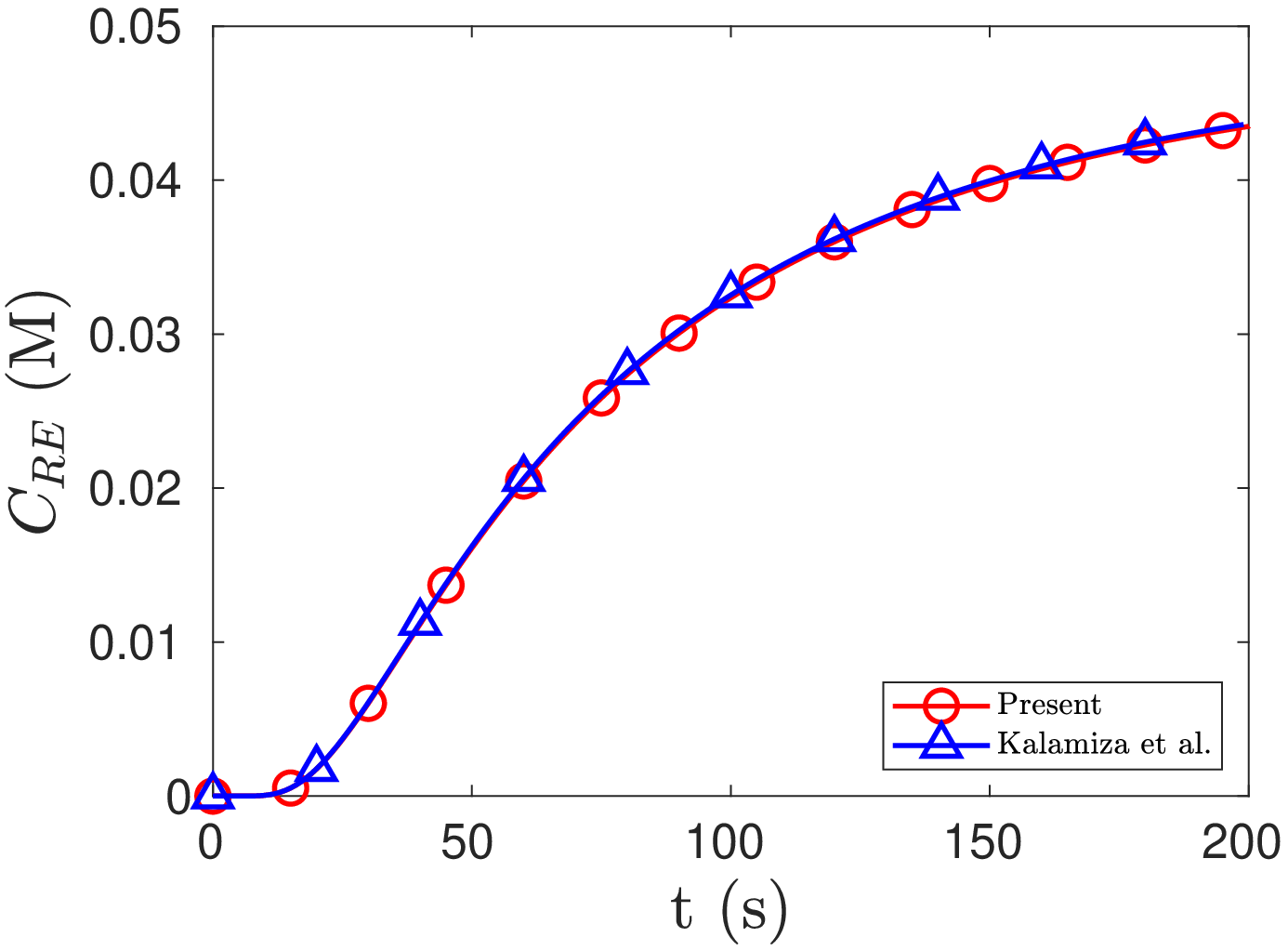}
		\\(a) & (b) \\		
	\end{tabular}
	\caption{The results for MTC and intracellular drug concentrations are compared with the result of the existing model \cite{Kalamiza20141950}. (a) To match MTC, $f_p=2.7\times 10^{-7}$ and $d_m=5\times 10^{-6}$ mm are taken. (b) The drug concentrations ($C_{RE}$) are calculated at (0.5, 0.5) for $P=0.0005$ mm s$^{-1}$, $D=0.001$ mm$^2$ s$^{-1}$ and $PN=1$.}
	\label{3comparison2}
\end{figure}
In order to validate the present model, a comparison is made between the present results and the results of the model proposed by Kalamiza et al. \cite{Kalamiza20141950}. The MTC that depends on membrane permeability is taken from the model \cite{Kalamiza20141950}. There is no pore resealing term present in the MTC of the Kalamiza's model.  For the sake of 
comparison with our results,  the pore resealing term $\left( \exp\left(-\frac{t}{\tau}\right)\right)$ is included in the MTC $\left(\frac{3Df_p}{d_mr_c} \right)$ used by kalamiza et al.  \cite{Kalamiza20141950} in their article,  which can be expressed as,
\begin{equation}\label{3mut_kalamiza}
\mu_k(t)=\frac{3Df_p}{d_mr_c}\cdot \exp\left(-\frac{t}{\tau}\right)=0.0032\cdot \exp\left(-\frac{t}{\tau}\right),
\end{equation}
where $f_p=2.7\times 10^{-7}$ and $d_m=5\times 10^{-6}$ mm.

Fig. \ref{3comparison2} presents the results for both MTCs given in the equations \eqref{3mu} and \eqref{3mut_kalamiza}.  It can be observed from the graphs that the results are in good agreement.
\subsubsection{Effects of drug permeability on drug transport}
\begin{figure}[h!]
	\centering
	\begin{subfigure}{0.5\textwidth}
		\centering
		\includegraphics[width=\linewidth]{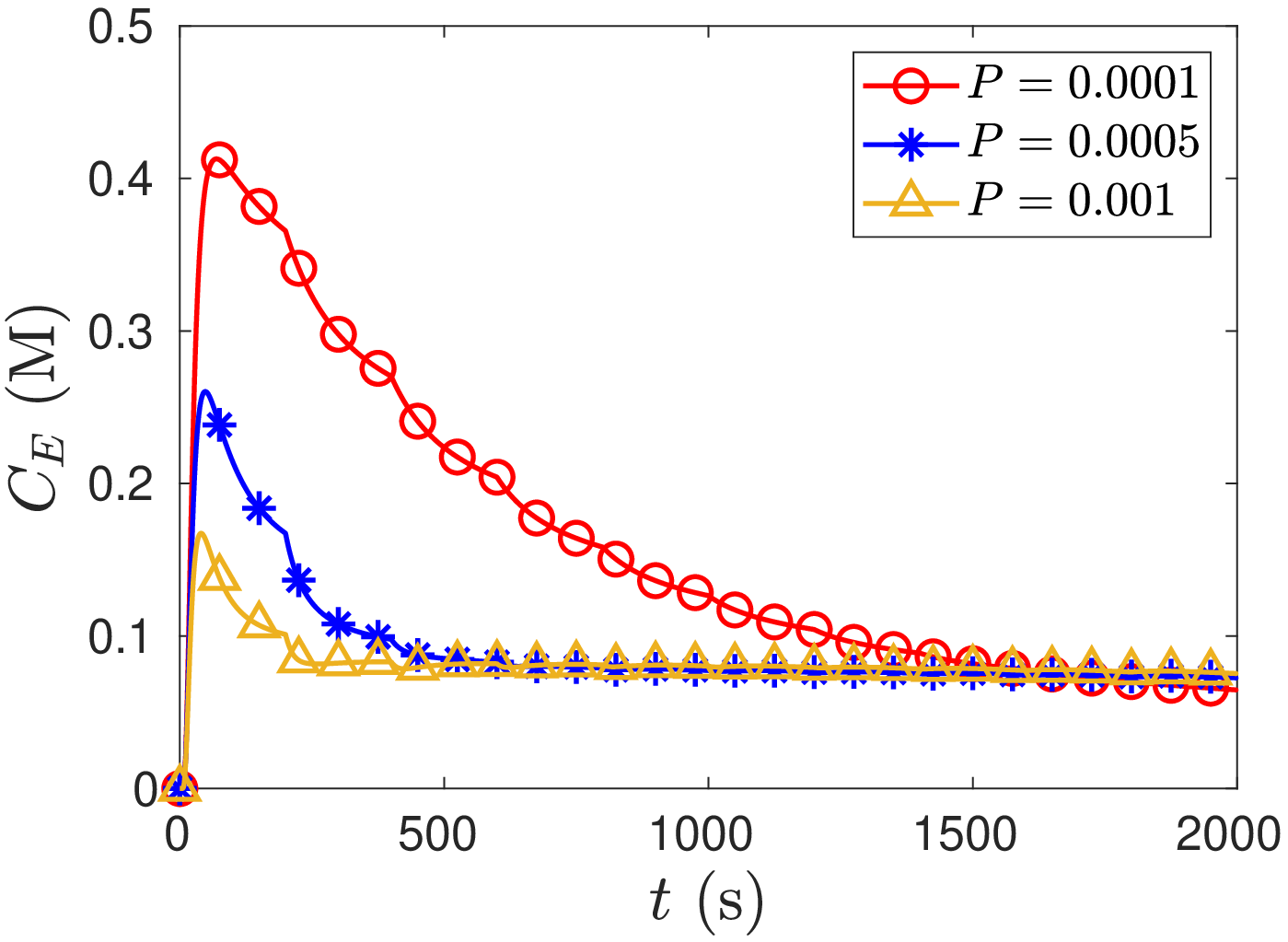}
		\caption{}
		\label{36sub1}
	\end{subfigure}%
	\begin{subfigure}{0.5\textwidth}
		\centering
		\includegraphics[width=\linewidth]{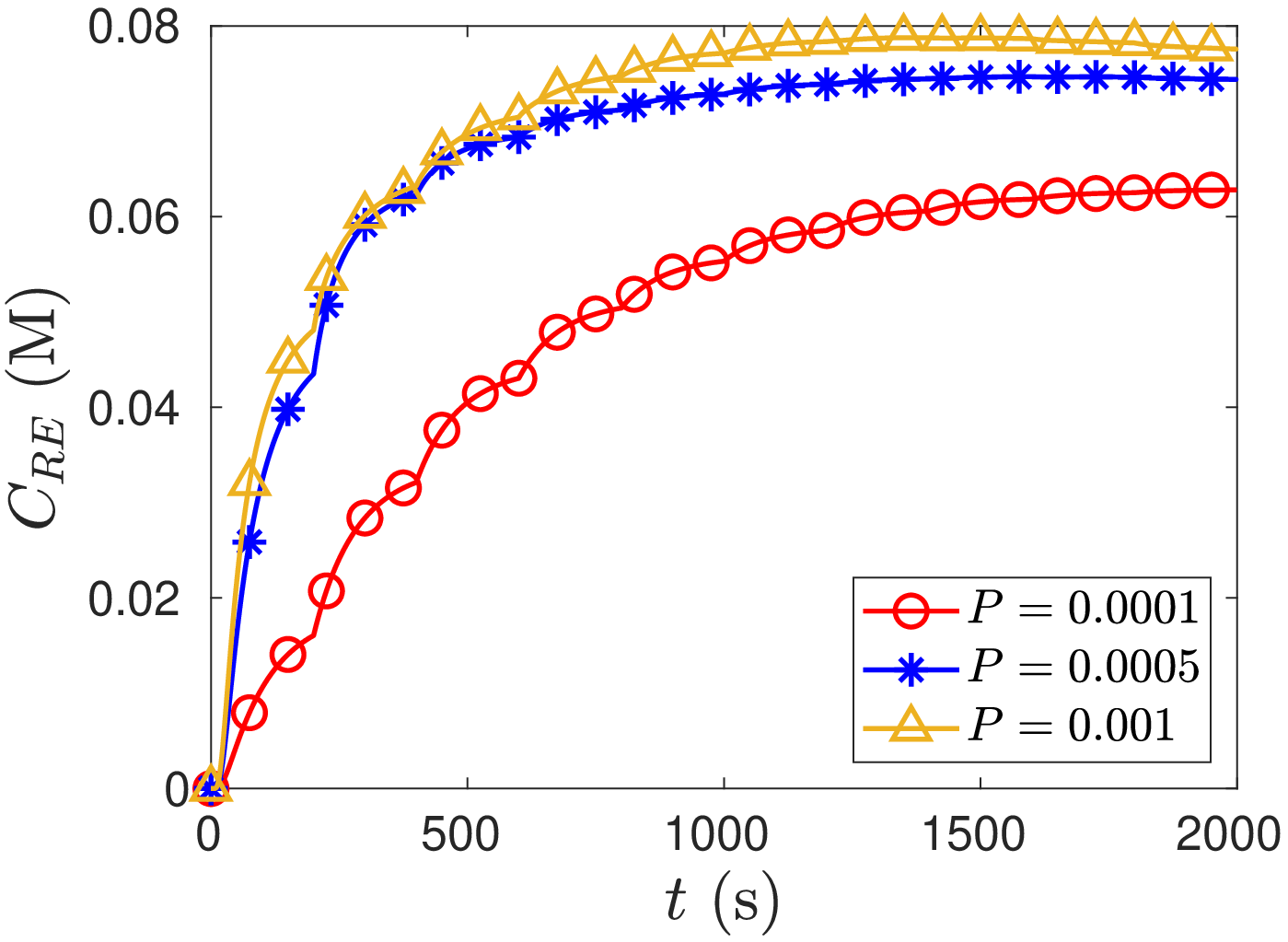}
		\caption{}
		\label{36sub2}
	\end{subfigure}
	\caption{Time-dependent drug concentration profiles for different permeability of drug.}
	\label{3CE_CRE_pchange_x5}
\end{figure}
The drug transport from extracellular to intracellular space depends on MTC, which is a function of drug permeability (see eq. \eqref{3mu}). So, the drug permeability parameter plays an important role in the drug uptake into the cells. The permeability dependency on drug concentrations ($C_E$, $C_{RE}$) is described through the Fig. \ref{3CE_CRE_pchange_x5}. It is noticed in the Fig. \ref{36sub1} that the drug concentration in the ECS decreases with time as there is a continuous drug uptake into the cells. Also, extracellular drug concentration decreases with the increase in drug permeability because intracellular drug uptake is more with the increase in permeability as shown in Fig. \ref{36sub2}.
\subsubsection{Effects of drug loss  on drug transport} 
\begin{figure}[h!]
	\centering
	\begin{tabular}{cc}
		\includegraphics[width=0.5\linewidth]{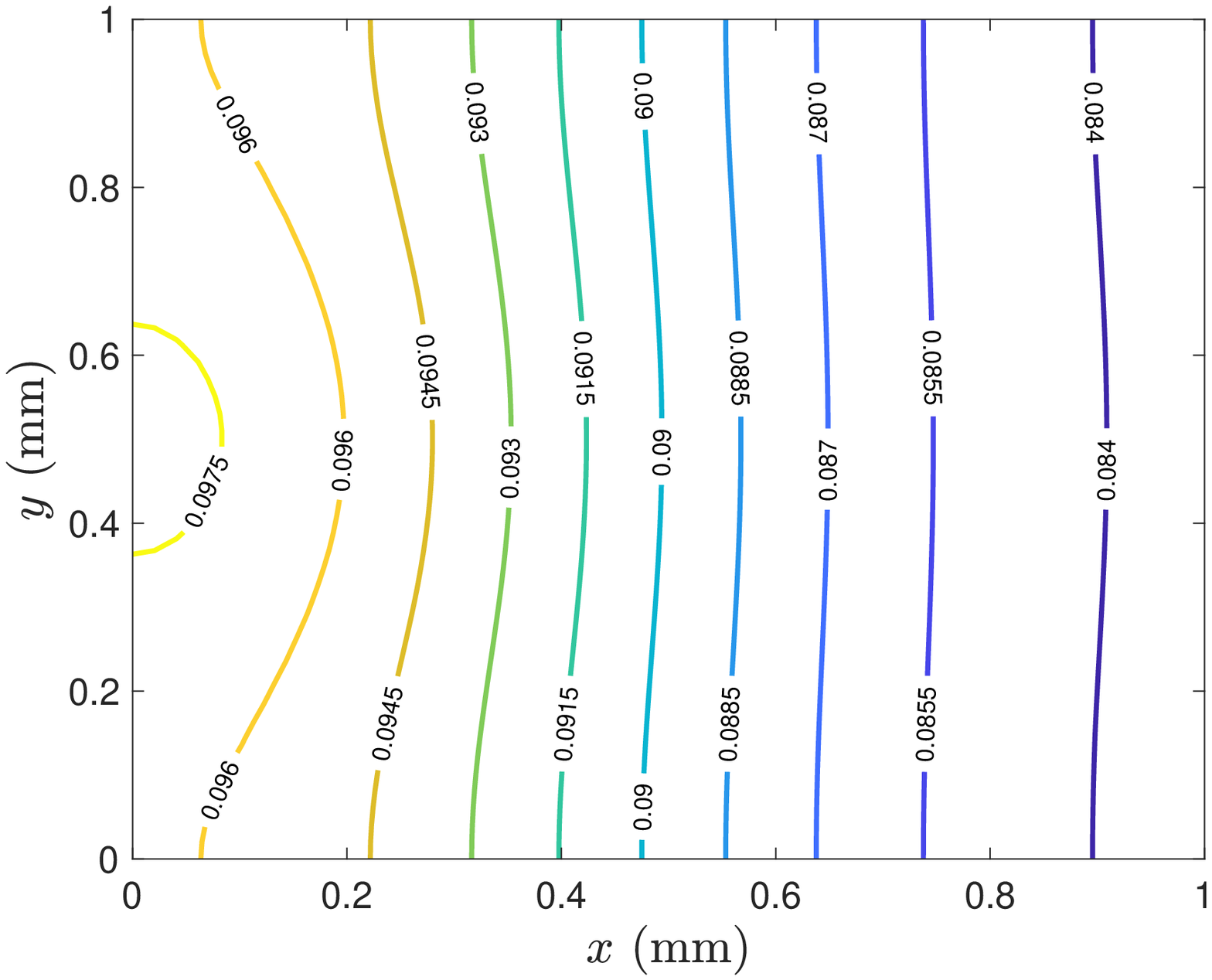}
		&
		\includegraphics[width=0.5\linewidth]{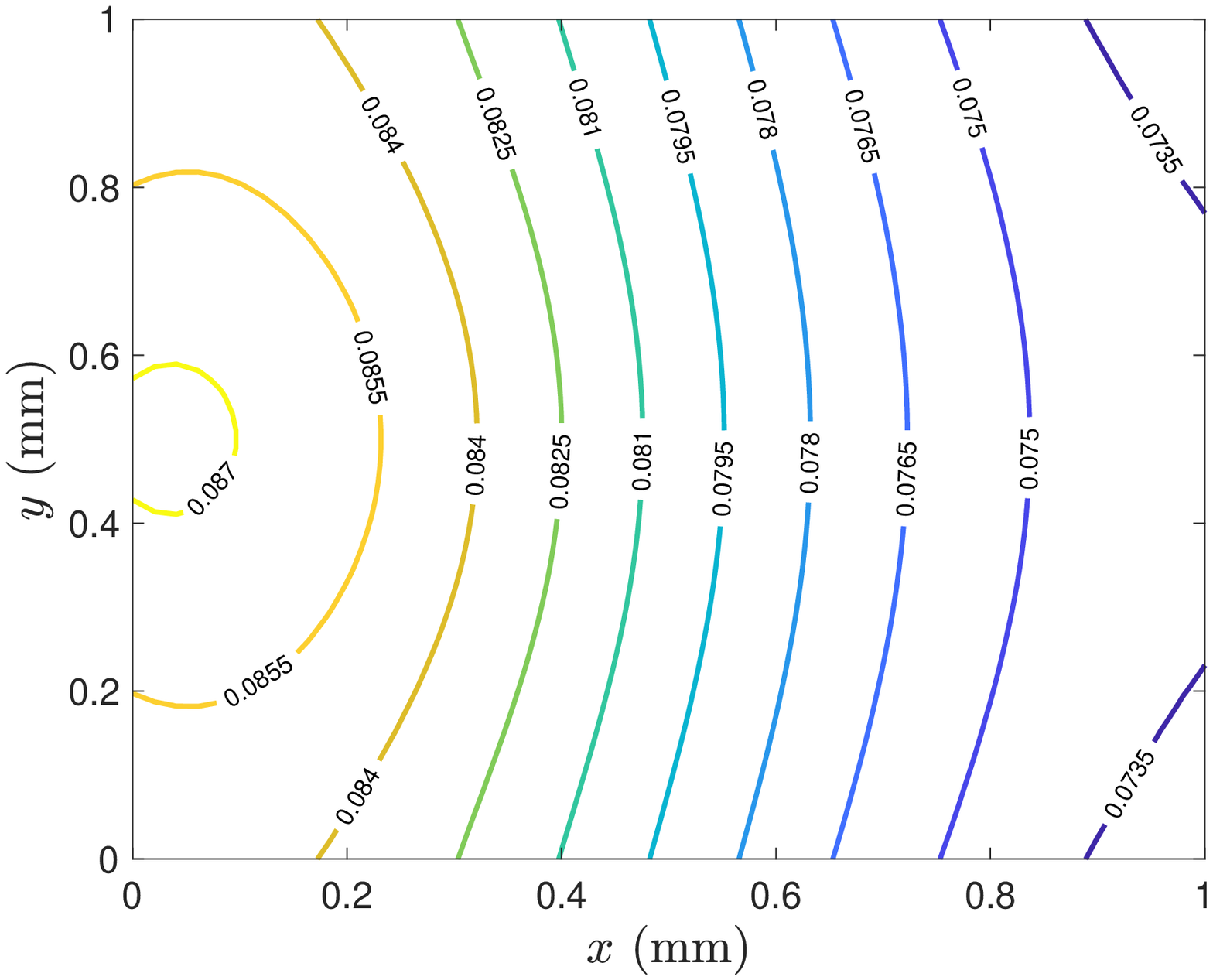}
		\\(a) $\beta=0$ & (b) $\beta=0.05$ \\	
		\includegraphics[width=0.5\linewidth]{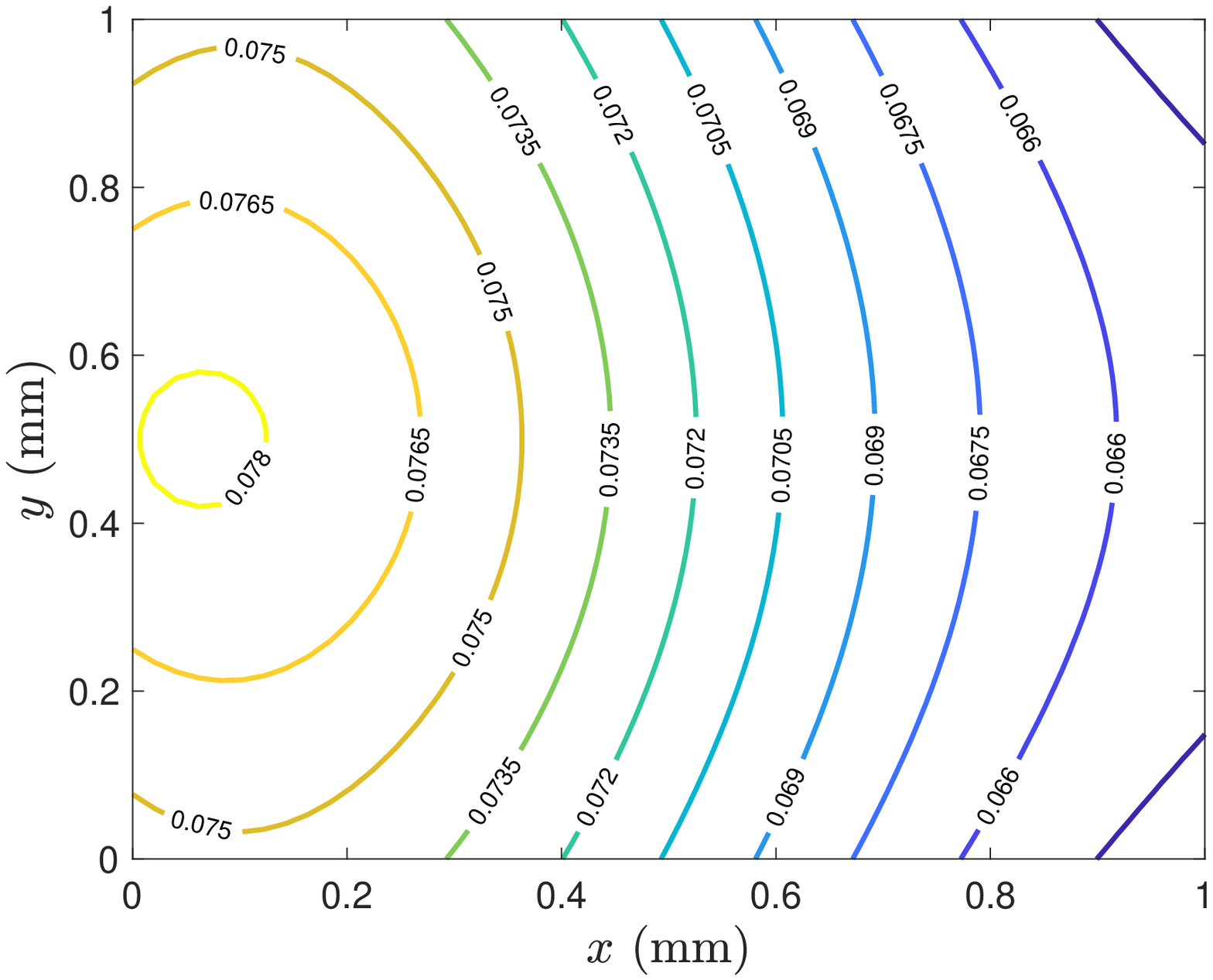}
		&
		\includegraphics[width=0.5\linewidth]{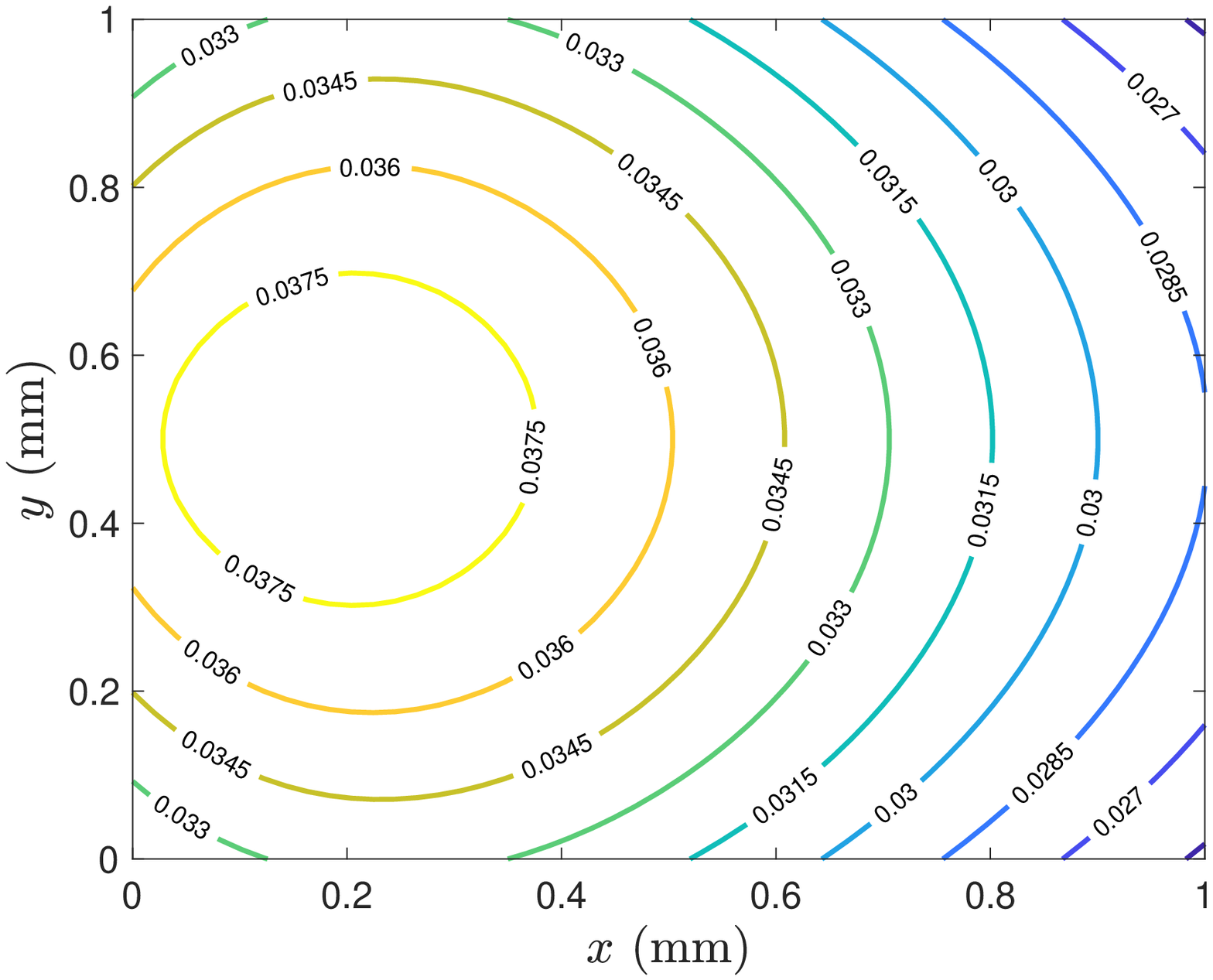}
		\\(c) $\beta=0.1$ & (d) $\beta=0.5$ 					
	\end{tabular}
	\caption{Drug distribution in the ECS after electroporation with 10 pulses of 1 ms for different values of $\beta$ (mm $^{-1}$).}
	\label{3CE_contour_beta_change}
\end{figure}
\begin{figure}[h!]
	\centering
	\begin{tabular}{cc}
		\includegraphics[width=0.5\linewidth]{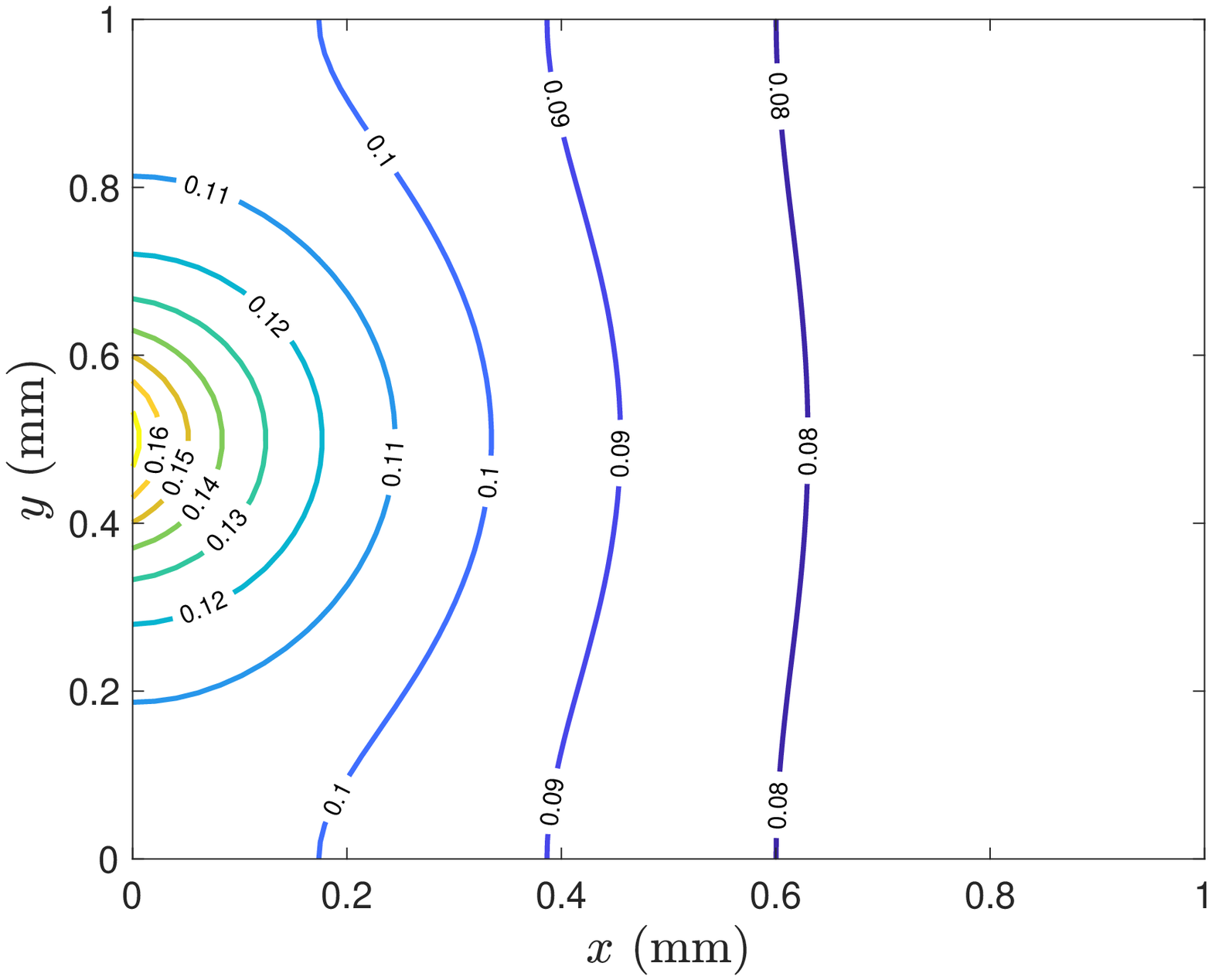}
		&
		\includegraphics[width=0.5\linewidth]{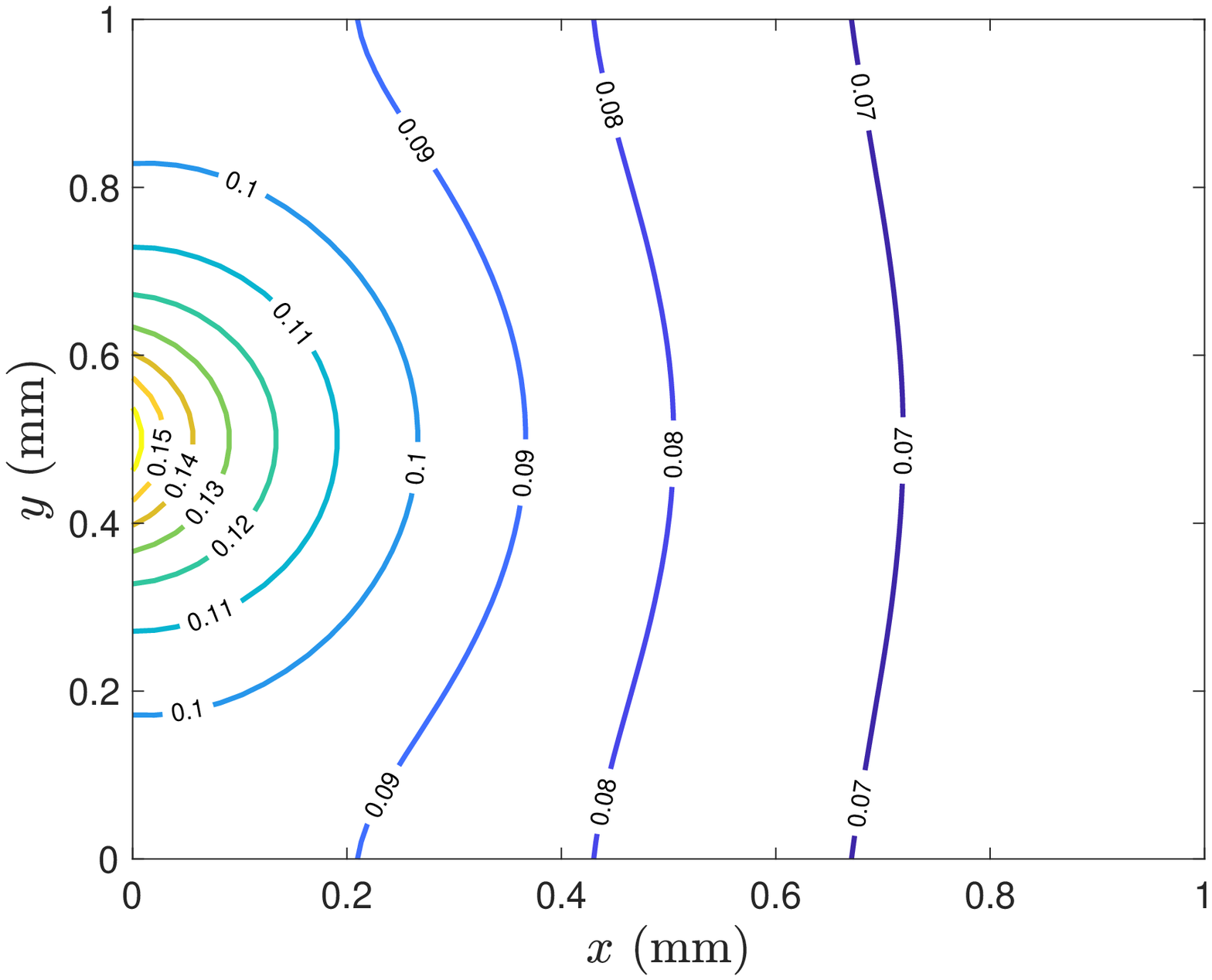}
		\\(a) $\beta=0$ & (b) $\beta=0.05$ \\	
		\includegraphics[width=0.5\linewidth]{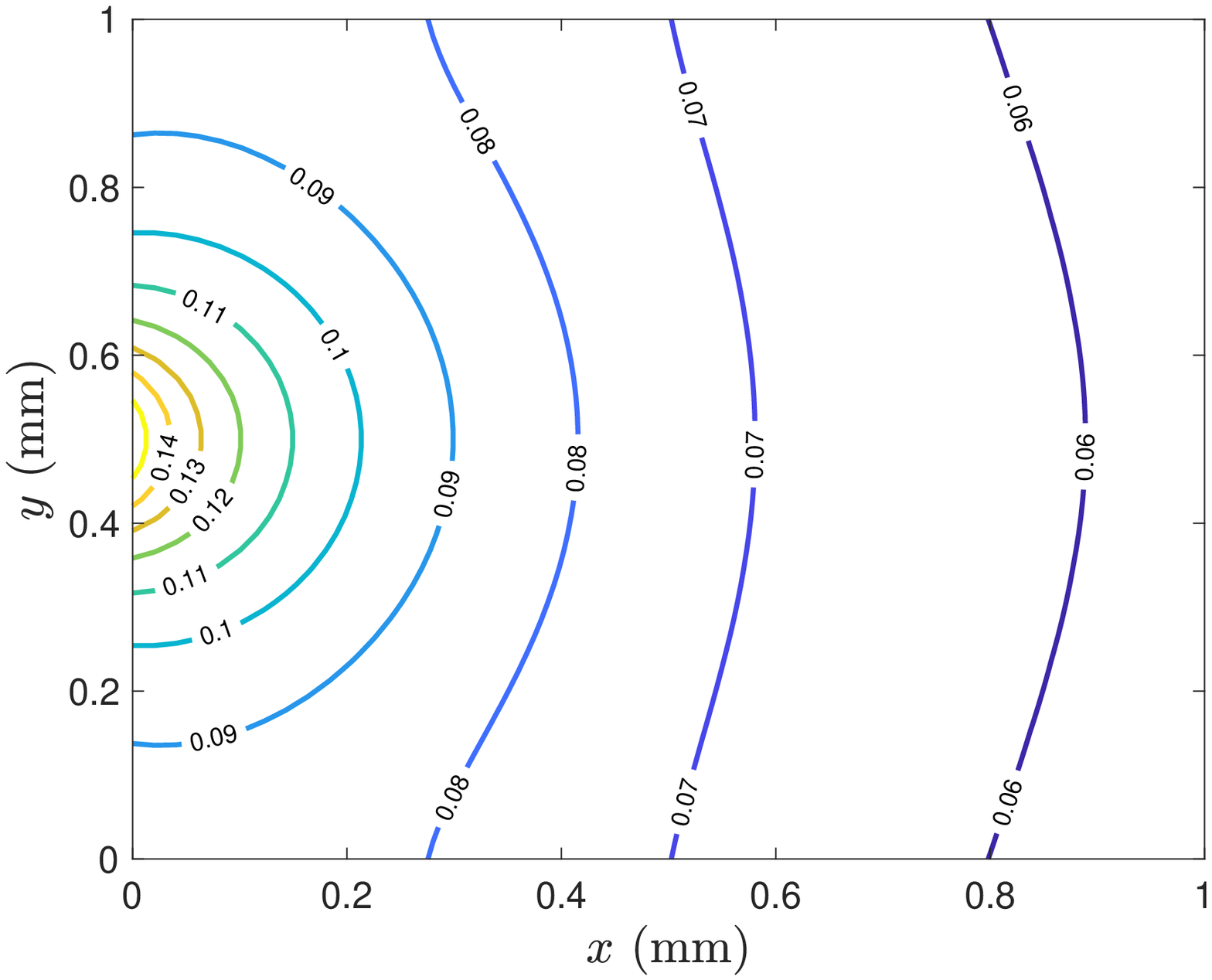}
		&
		\includegraphics[width=0.5\linewidth]{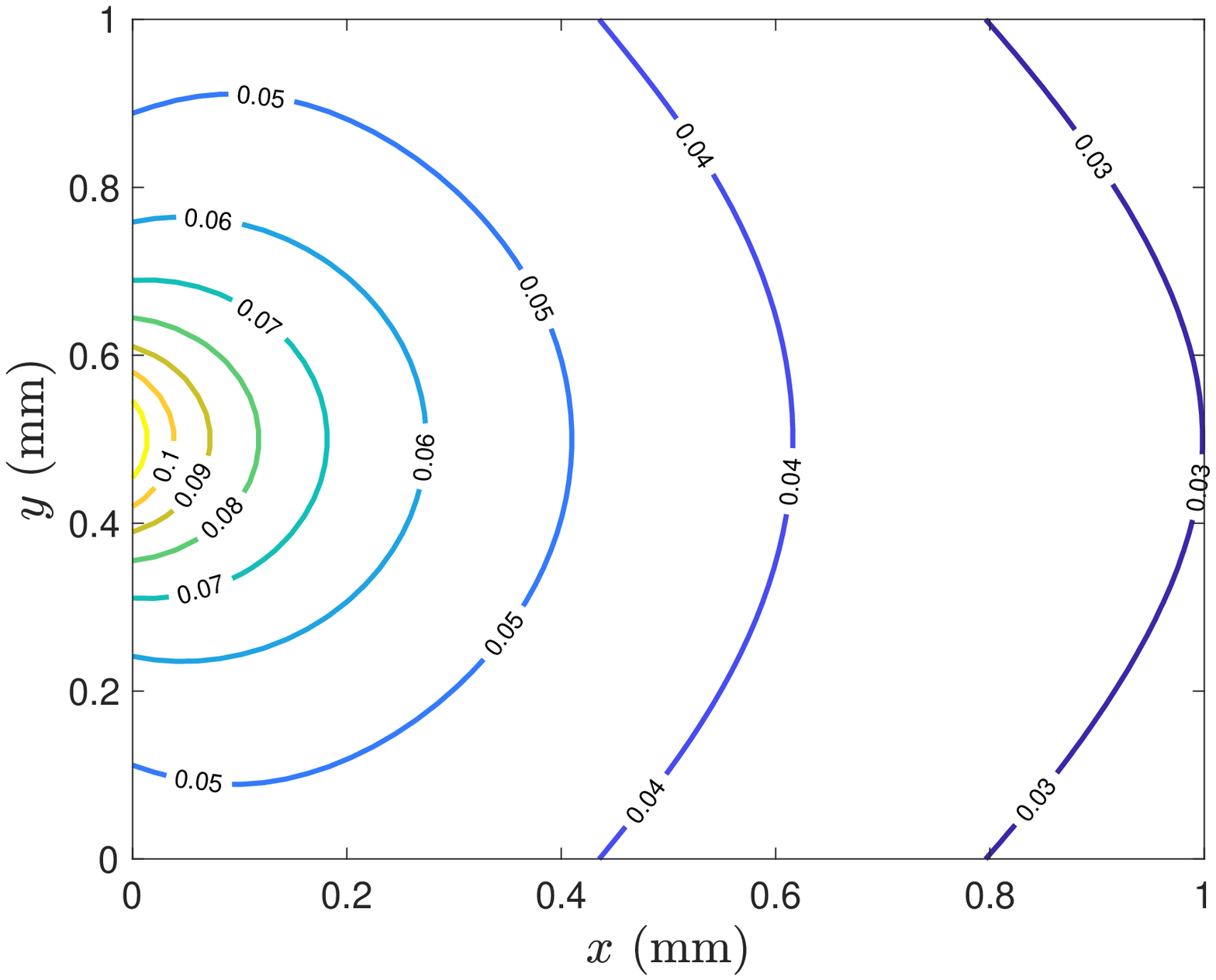}
		\\(c) $\beta=0.1$ & (d) $\beta=0.5$ 					
	\end{tabular}
	\caption{Drug distribution in the ICS after electroporation with 10 pulses of 1 ms for different values of $\beta$ (mm$^{-1}$).}
	\label{3CRE_contour_beta_change}
\end{figure}
\begin{figure}[h!]
	\centering
	\begin{tabular}{cc}
		\includegraphics[width=0.5\linewidth]{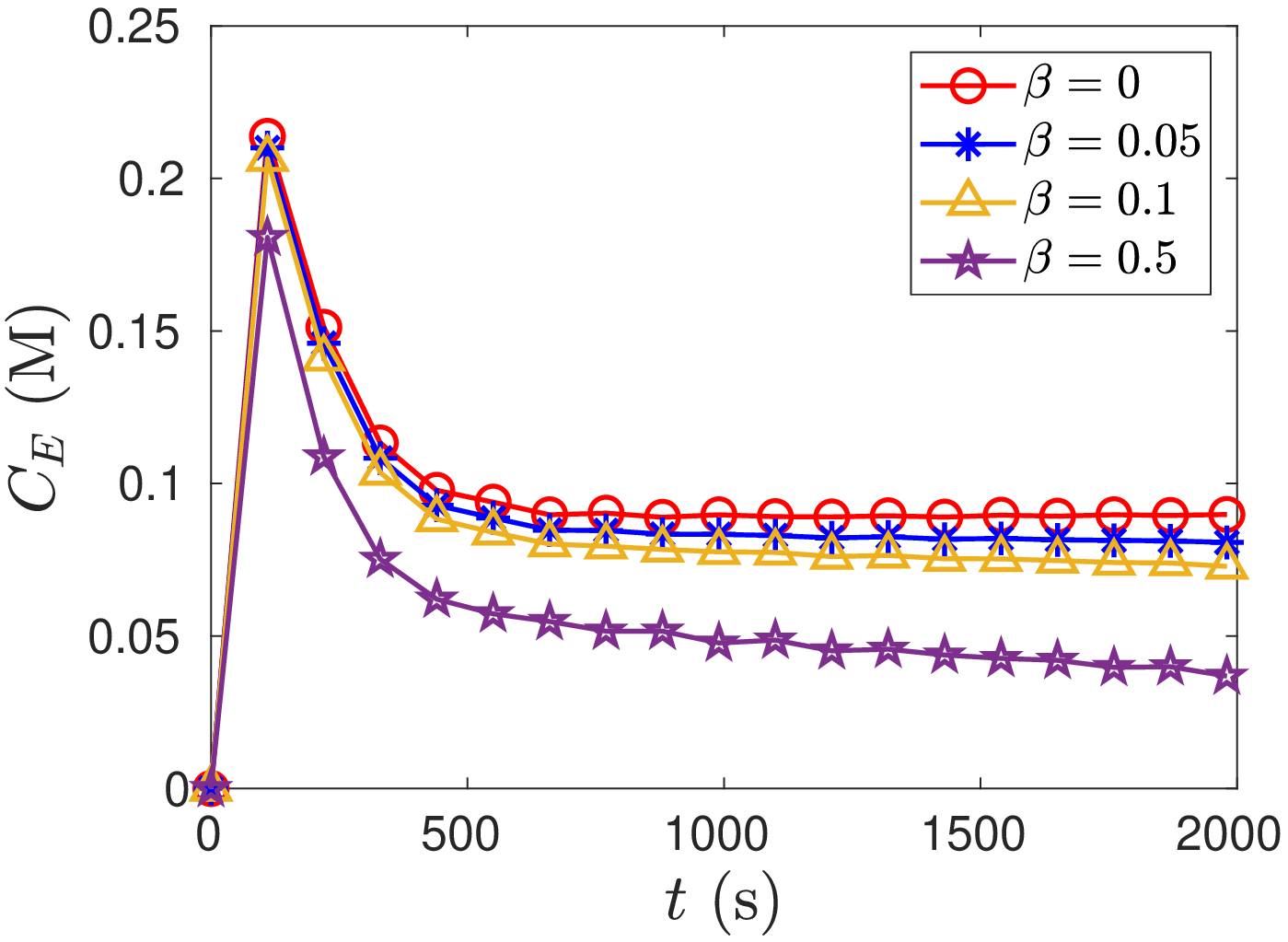}
		&
		\includegraphics[width=0.5\linewidth]{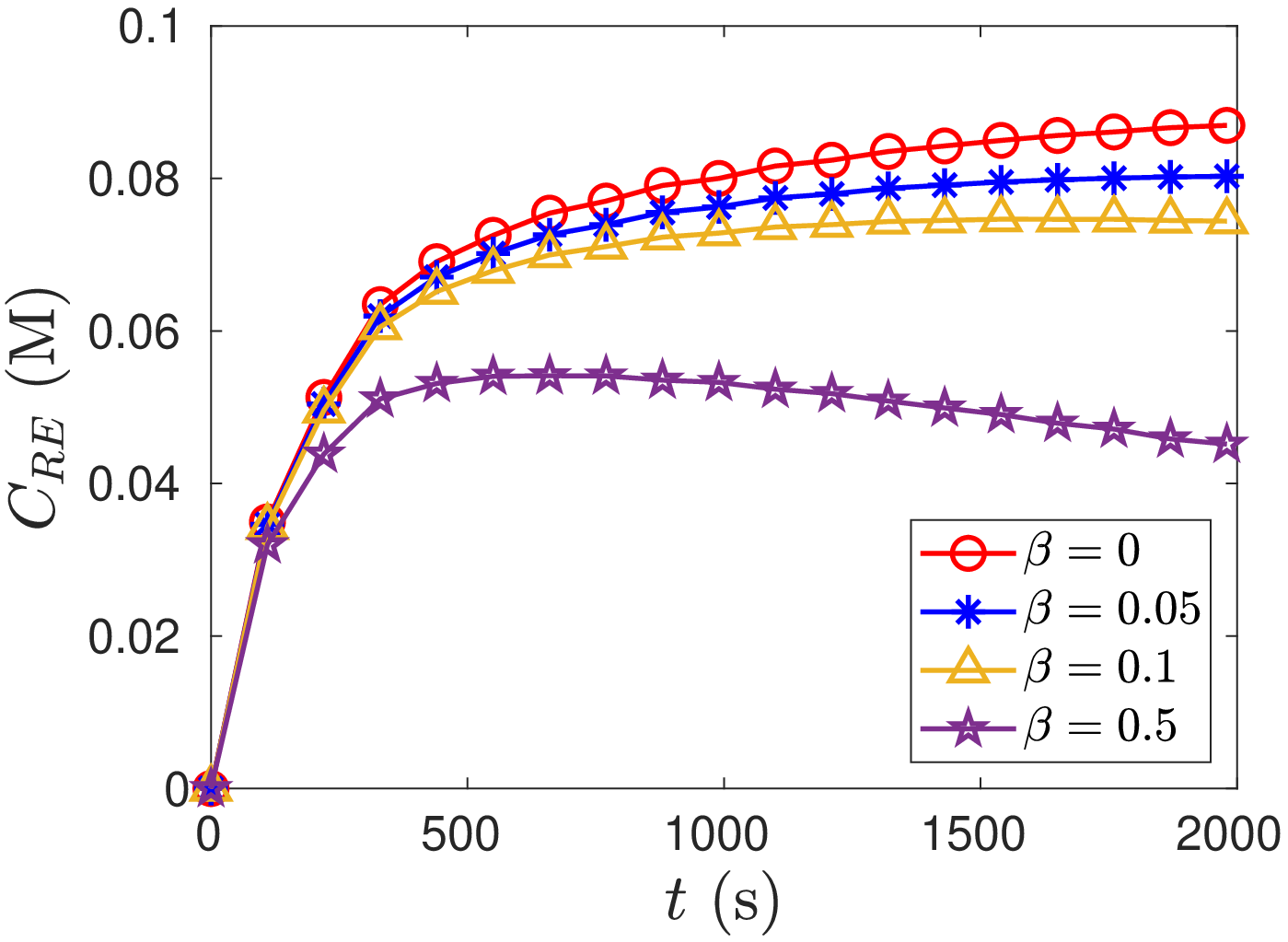}
		\\(a) & (b) \\	
	\end{tabular}
	\caption{Time-dependent concentration profiles in (a) ECS and (b) ICS for different values of $\beta$ (mm$^{-1}$).}
	\label{3CE_CRE_betachange}
\end{figure}
The drug penetration throughout the extracellular media in the tissue for different rates of drug loss ($\beta$) at the tissue boundaries is described in Fig. \ref{3CE_contour_beta_change}. The results are plotted after performing 10 pulses of 1 ms with a rest time  200 s between two pulses. It is clearly noticed that the drug penetration in the tissue is affected with the increase of $\beta$. This is obvious as more drugs move outside from the tissue domain with the increase in $\beta$. It is also observed in the Fig. \ref{3CE_contour_beta_change}d that the amount of drugs reach into the ECS is less for high value of $\beta$. However, higher amount of drugs enters into the ECS for lower value of $\beta$ and maximum amount of drugs reaches for $\beta=0$ (i.e., no drug loss from tissue boundaries)  which is shown in the Fig. \ref{3CE_contour_beta_change}a.

Fig. \ref{3CRE_contour_beta_change} shows the spatial changes of cellular drug uptake in the tissue for different values of $\beta$. It is observed from the graphs that the drug uptake into the cells decreases with the increase in $\beta$. It may be due to the fact that fewer drugs reach in the ECS around cells for higher values of $\beta$, as observed in \ref{3CE_contour_beta_change} so that less amount of drugs enters into the cells from ECS. Also, it is noticed from the contour plots that the amount of drugs intake into the cells decreases along x-axis from left to right as the drug source is located at the left boundary.

The temporal changes of drug concentrations in ECS and ICS for different values of the parameter $\beta$ are illustrated by the Fig. \ref{3CE_CRE_betachange}. According to Fig. \ref{3CE_CRE_betachange}a, the extracellular drug concentration increases rapidly within the first pulse as a result of rapid drug diffusion in the ECS. It starts decreasing when the drug enters the cells and continues to do so until it reaches an equilibrium state. Fig. \ref{3CE_CRE_betachange}b shows that the intracellular drug profiles  increase throughout the time for all values of $\beta$ except $\beta=0.5$. For the less values of $\beta$ ($=0, 0.05, 0.1$), the concentration increases due to continuous drug intake into the cells with the application of repeated pulses and low rate of drug losses.  For the higher value of $\beta$ ($=0.5$), the cellular drug uptake increases initially but after a certain time ($\approx 500$ s), it decreases due to huge drug losses.

\subsubsection{Effects of pulse number on drug transport}
\begin{figure}[h!]
	\centering
	\begin{tabular}{cc}
		\includegraphics[width=0.5\linewidth]{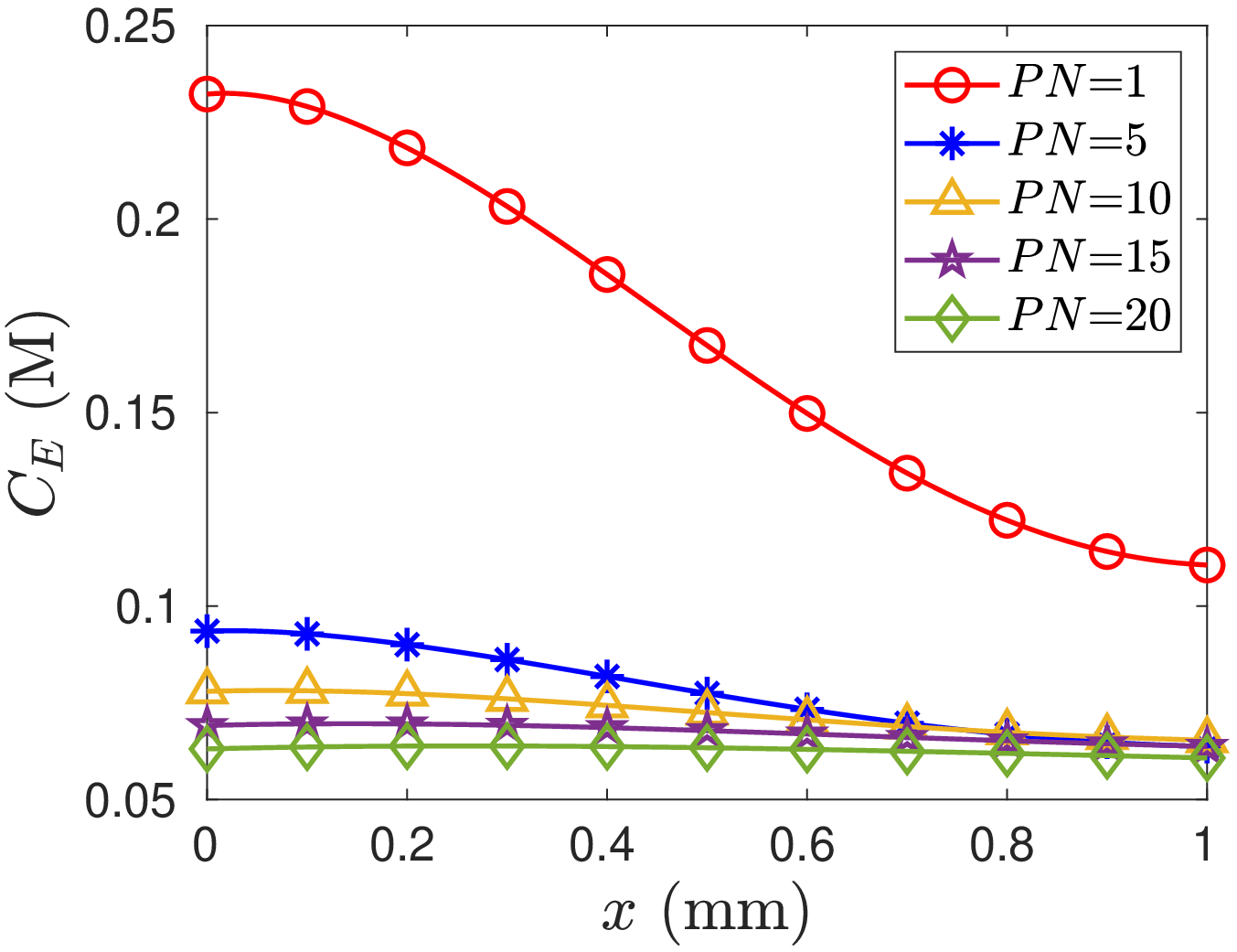}
		&
		\includegraphics[width=0.5\linewidth]{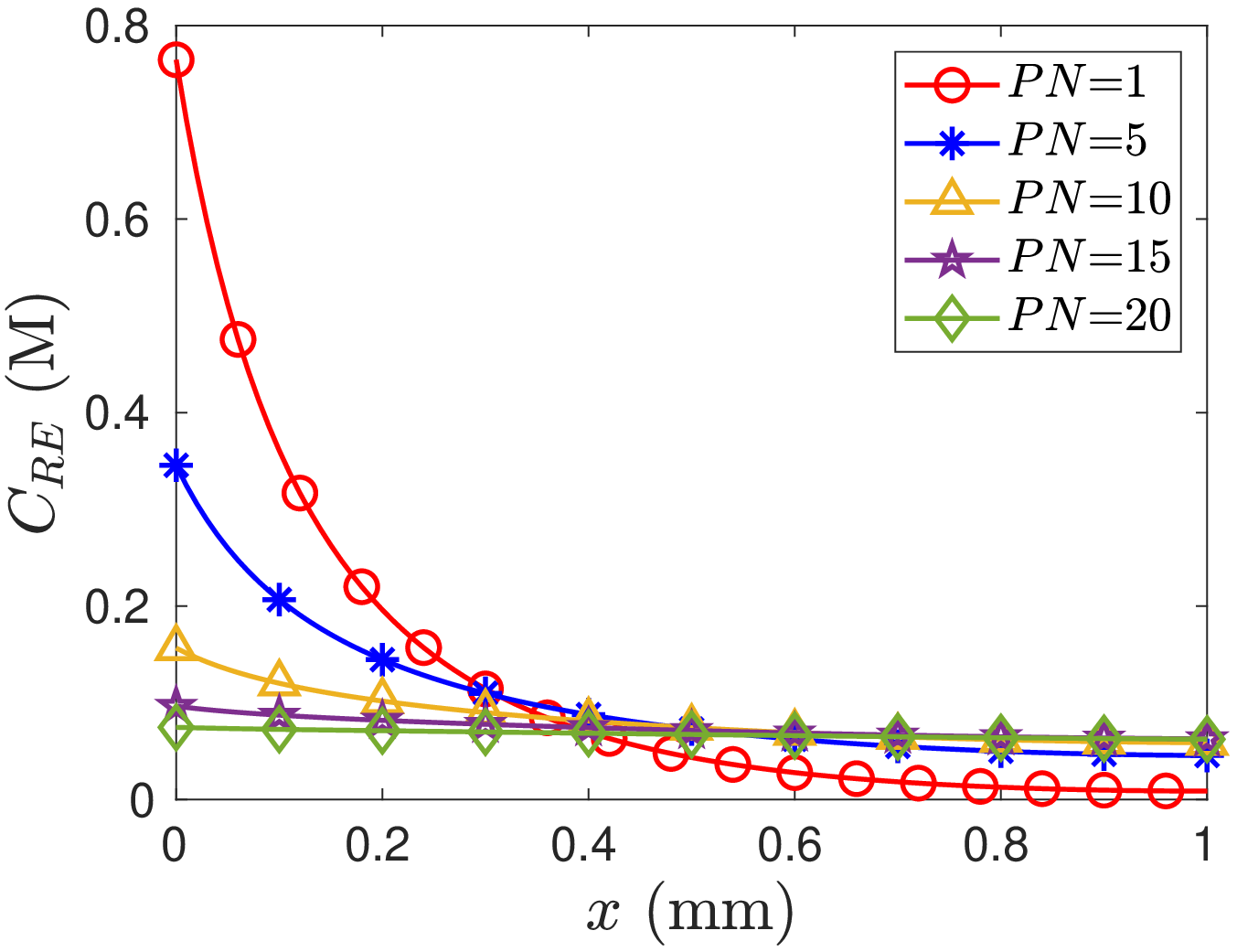}
		\\(a) & (b) \\	
	\end{tabular}
	\caption{Concentration changes along x-axis in (a) ECS and (b) ICS at $y=0.5$ for $\beta=0.1$ mm$^{-1}$ and different pulse numbers.}
	\label{3CE_CRE_PN_change_beta_01}
\end{figure}

\begin{figure}[h!]
	\centering
	\includegraphics[width=0.9\linewidth]{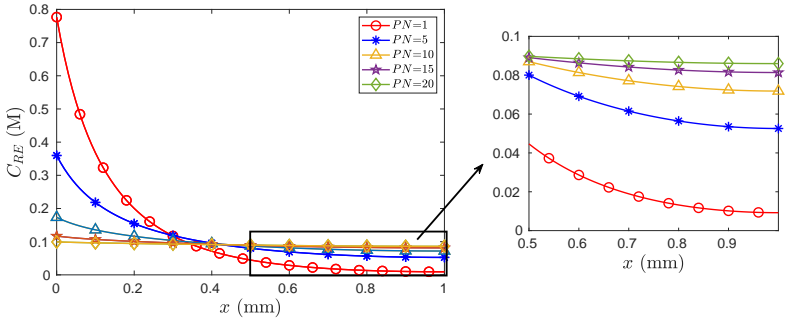}
	\caption{Concentration changes along x-axis in ICS  at $y=0.5$ for $\beta=0$ per mm and different pulse numbers.}
	\label{3CE_CRE_PN_change_beta_0}
\end{figure}
In this model, multiple pulses are applied to enhance the drug transport rate into the cells. If a single pulse is applied and mass transport occurs for a long duration, the mass transfer rate slows down due to pore resealing effects. Fig. \ref{3CE_CRE_PN_change_beta_01} describes the changes of drug concentrations in ECS and ICS with the variation of the pulse number. It is seen from Fig. \ref{3CE_CRE_PN_change_beta_01}a that the extracellular concentration decreases with the increase in pulse number. The possible reasons are: (i) the increment of drug intake into the cells with the increase in pulse number; (ii) the drug loss from the tissue boundaries (as $\beta=0.1$).
In Fig. \ref{3CE_CRE_PN_change_beta_01}b, it is noticed that the intracellular drug concentration decreases with pulse number in the region $0\le x<0.4$ and increases in the region $0.4\le x \le 1$ of the tissue.  This opposite behavior before and after $x=0.4$ is noticed due to continuous drug diffusion from left to right and ongoing drug introduction into the cells  in the domain $0.4\le x \le 1$.  Fig. \ref{3CE_CRE_PN_change_beta_0} is plotted when there is no drug loss (i.e, $\beta=0$) to show the actual effects of pulse number. One observation is made from the figure that multiple pulses are required for the uniform drug distribution into the targeted cells of the tissue. It is concluded from the above discussion that the drug concentrations in both the regions ECS and ICS reach a stable state after the application of a certain number of pulses ($PN=20$, in this case). Hence, for this model, the maximum number of pulses is 20, which is enough to inject a sufficient quantity of drug into the cells uniformly.


\section{Conclusion} 
The present study proposes an improved mathematical model for drug delivery in which a point source drug is introduced at the tissue boundary. The effects of drug loss from the boundary of the targeted tissue is also emphasized through this model. The following significant observations are made from this study:
\begin{itemize}
	\item The prescribed model is capable of describing a complete overview of drug transport phenomena in the targeted tissue with reversible electroporation. 
	
	\item In the electroporation based drug delivery, the Dirac-delta function is used  to represent the point source drug more appropriately.
	
	\item The model deals with Neumann boundary condition to incorporate the drug elimination process from the tissue.
	
	\item It is found that the drug loss from the tissue boundaries is an important factor in drug transport.  Through this model, the mechanism of drug dose is analyzed for the first time.
	
	\item It is observed that the drug uptake into the cells increases with applications of multiple pulses.
	
	\item It is also seen that the cellular drug uptake is improved for smaller values of $\beta$. 
	
	\item The drug uptake into the cell becomes slow, even if multiple pulses are applied, due to huge drug loss at the tissue boundaries.
	
	\item Uniform drug distribution into the cells is done successfully. So, the model configuration is a successful tool to treat the whole targeted tissue by introducing drug as a medicine.
\end{itemize}
The proposed model could be helpful in pharmaceutical industry and medical sciences.

%
%
%
%

\bmhead{Acknowledgments}

The first author (Nilay Mondal) acknowledges the support received from the Indian Institute of Technology Guwahati, Assam, India.

\bibliography{sn-article.bbl}


\end{document}